%
\documentstyle{amsppt}
\voffset=-50pt
\magnification=1200
\catcode`\@=11
\def\nologo{\let\logo@\empty}
\catcode`\@=\active
\nologo
\def\lf{\left}
\def\ri{\right}
\def\H{\Bbb H}
\def\R{\Bbb R}
\def\D{\Bbb D}
\def\C{\Bbb C}
\def\Z{\Bbb Z}
\def\i{\sqrt{-1}}
\def\dist{\text{dist}}
\def\e{\exp}
\def\z{\zeta}

\def\ol{\overline}
\def\bqd{\text {BQD}(\H^2)}
\def\t{\Cal T}
\def\cS{{\Cal S}}
\def\b{\frak B}
\def\tm{\tilde\mu}
\def\tn{\tilde\nu}
\def\tf{\tilde F}
\topmatter
\title
 Hopf differentials and   the images of harmonic maps   \endtitle
\rightheadtext{ }
\leftheadtext{ }
\author
Thomas K. K.  Au,  Luen-Fai Tam,\footnotemark \, 
and Tom Y. H. Wan\footnotemark
\endauthor
\footnotetext"$^1$"{Research was partially supported by a grant from the Research Grants Council of the Hong Kong Special Administrative Region (Project No. CUHK4217/99P)}
\footnotetext"$^2$"{Research was partially supported by a grant from the Research Grants Council of the Hong Kong Special Administrative Region (Project No. CUHK4199/97P)}
\affil
 Department of Mathematics\\
The Chinese University of Hong Kong\\
Shatin, Hong Kong 
\\
thomasau\@math.cuhk.edu.hk\\
lftam\@math.cuhk.edu.hk\\
tomwan\@math.cuhk.edu.hk
\\
\endaffil
\date
Revised Oct, 2000
\enddate
\endtopmatter

In \cite{Hz}, Heinz proved that there is no harmonic diffeomorphism from the unit disk $\D$ onto the complex plane $\C$. The result was   generalized by Schoen  \cite{S} and  he proved  that there is no harmonic diffeomorphism from the unit disk onto a complete surface of nonnegative curvature. Unlike conformal or quasi-conformal maps between Riemann surfaces, the inverse of a harmonic map is not harmonic in general. Hence it is an interesting question whether there is any harmonic diffeomorphism from $\C$ onto $\D$ equipped with the Poincar\'e metric. In fact a general form of this  question was formulated by Schoen  \cite{S} as follows: ~{\sl Is it true that Riemann surfaces which are related by a harmonic diffeomorphism are necessarily quasi-conformally related?} 

Let us first recall some facts on harmonic maps between surfaces. Let $\Sigma_1$  and
$\Sigma_2$ be two Riemann surfaces with conformal metrics
$\rho^2(z)|dz|^2$ and
$\sigma^2(h) |dh|^2$ respectively. The harmonic map equation for maps
from
$\Sigma_1$ into
$\Sigma_2$ can be written as 
$$
h_{z\ol z}+2(\log \sigma)_hh_zh_{\ol z}=0.
$$
Define $||\partial h||= \rho^{-1}\sigma |h_z|$, and $||\ol\partial h|| =\rho^{-1}\sigma |h_{\ol z}|$. Hence $||\partial h||$ and $||\ol\partial h||$ are the norms of the $(1,0)$-part and $(0,1)$-part of $dh$. The energy
density of $h$ is given by $e(h)=||\partial h||^2+||\ol\partial h||^2$,
and the Jacobian of $h$ is given by $J(h)=||\partial h||^2-||\ol\partial
h||^2$.   The Hopf
differential of $h$ is defined as $\phi dz^2=\sigma^2(h)h_z\ol h_z dz^2$, which is the $(2,0)$-part of $h^*\lf(\sigma^2(h)|dh|^2\ri)$.
It is well known that if $h$ is harmonic then $\phi dz^2$ is a holomorphic
quadratic differential defined on $\Sigma_1$, see \cite{C-G}.  If $h$ is an orientation preserving local  diffeomorphism, then  $J(u)>0$, which
implies that
$||\partial h||>0$ everywhere, and that
$$
e^{2w} >|\phi| $$ 
where $w=\log||\partial h||$.

In \cite{Wn}, the third author proved that orientation preserving harmonic diffeomorphisms on the  hyperbolic plane $\H^2$  can be parametrized by their Hopf differentials, provided that they can be realized as the Gauss maps of constant mean curvature cuts in the Minkowski space $M^{2,1}$. The last condition is equivalent to the fact that $|\partial h |^2|dz|^2$ is a complete metric on $\D$, where $|\partial h|$ is the norm of $\partial h$ with respect to the Euclidean metric on the domain and the Poincar\'e metric on the target.  The result was generalized to harmonic maps from $\C$ into $\H^2$ in  \cite{W-A} and to more general  surfaces in \cite{T-W 1}.   

Hence in order to study the behaviors of harmonic maps from $\C$ or $\D$ into $\H^2$ it is useful to study their Hopf differentials. In \cite{Hn}, Han proved that if $h$ is a harmonic diffeomorphism from $\C$ into $\H^2$ whose Hopf differential is a polynomial, then the closure of $h(\C)$ in $\overline{\H^2}=\H^2\cup\partial\H^2$ is the convex hull of a totally disconnected closed set in $\partial \H^2$, provided that $|\partial h |^2|dz|^2$ is   complete on $\C$. In particular, $h$ is not surjective. Here $\partial \H^2$ is the geometric boundary of $\H^2$. Later in \cite{HTTW}, it was proved that the closure of the image of an orientation preserving  harmonic diffeomorphism $h$ from $\C$ into $\H^2$ is an ideal polygon with exactly $m+2$ vertices on $\partial \H^2$ if and only if the Hopf differential $\phi~dz^2$ is a polynomial of degree $m$, i.e. $\phi$ is a polynomial of degree $m$. Note that by \cite{Wn, T-W 1}, we know that $\phi$ is of degree no greater than $m$ if and only if $h$ is of polynomial growth of  degree at most $m/2+1$.

In higher dimension, one cannot expect that such a clean statement continues to hold. However, in \cite{L-W 1}, Li and Wang were able to generalized   part of the above result for a much more general class of manifolds. They proved that if $M^n$ is a complete manifold with nonnegative Ricci curvature and $N$ is a Cartan-Hadamard manifold with sectional curvature pinched between two negative constants, then the closure in $N\cup\partial N$ of  the image of  harmonic map from $M^n$ into $N$ with polynomial growth of degree at most $\ell$  is in the convex hull of finitely many points on the geometric boundary $\partial N$ of $N$. Moreover, the number of points is bounded by a constant depending only on $n$ an $\ell$. Actually, they only assumed that $M$ satisfies the so-called weak volume growth  condition and weak Poincar\'e inequality. In \cite{L-W 3}, they also obtained a sharp estimate for the number of points on the ideal boundary in case $M$ is a complete surface with finite total curvature.

All these results in \cite{Hn, HTTW, L-W 1, 3} are  under the assumption that the harmonic map is of polynomial growth. In this paper, we want to study harmonic maps from $\C$ into $\H^2$ which grow faster than polynomial. We will study the images of the harmonic maps by a careful study of their Hopf differentials.

 First we prove that if $h$ is  an orientation preserving harmonic diffeomorphism with Hopf differential $\phi~dz^2$ such that  $\phi$ is of one of the following forms then $h$ is not surjective:
\roster
\item $\phi=P_1\e\left[P_2\e\left[\cdots \e\left[P_k\e(Q)\right]\cdots\right]\right]$, where $P_j$ and $Q$ are polynomials (Theorem 1.2);
\item $\phi=P\e(Q) $ where  $Q(z)=z^n+\sum_{j=1}^{n}a_jz^{n-j}$ is a polynomial of degree $n\ge 1$,   $P$ is entire with order $\rho<n$ and   
$$
\Sigma\cap\{z|\ |z|>R_0 \ \text{and}\ -\delta<\arg z <\delta\}=\emptyset,
$$
  for some  $\frac\pi n>\delta>0$ and $R_0>0$, where  $\Sigma$ is the set of all zeros of $P$ (Theorem 1.1);
\item $\phi=(f')^2$ where $f$ is entire with  no finite asymptotic value in the domain
$$
\Cal R=\lf\{z|\ \frac\pi2-\delta<\arg z<\frac\pi2+\delta,\ \text{and}\ | z|>R\ri\} $$
for some  $\delta>0$ and $R>0$ and  $f'(z)\neq0$ for all $z$ in $f^{-1}(\Cal R)$ (Theorem 1.3).
\endroster
In (2), $\phi$ is of finite order, in (1) $\phi$ is of infinite order and there is no growth condition in (3). Note that if $\phi$ is of finite order then $\phi=P\e(Q)$ with $Q$ being a polynomial and $P$ is entire.

As mentioned above, if $|\partial h|^2|dz|^2$ is complete on $\C$ and if $\phi$ is a polynomial, then   the image of $h$ is an ideal polygon with finitely many vertices  at $\partial \H^2$  \cite{HTTW}. If $\phi$ is not a polynomial, then $\overline{h(\C)}\cap \partial\H^2$ must consist of infinitely many points by \cite{HTTW} again and in this case, the image set is much more difficult to be described. In the second part of this paper, we want to describe the images of harmonic maps under similar assumptions as in (1) or (2) above. We prove that if $\phi=P\e(Q)$ where $P$ and $Q$ are polynomials, then $\overline{h(\C)}\cap \partial\H^2$ is countable and consists of exactly $n$ accumulation points, where $n=\deg Q$ (Theorem 3.1). In fact, one can relax the condition that $P$ is polynomial. If we assume that $P$ is entire with order less than $n$ and the zeros of $P$ are well distributed, then the same conclusion holds. 
Next we consider the case that  $\phi(z)=P(e^z)$ where $ P(t)=\sum_{k=-m}^{n} a_kt^k$ and is non-constant. In this case $\phi$ is of order one. It is interesting to know that  under this assumption on $\phi$,  $\overline{h(\C)}\cap \partial\H^2$ has only one accumulation points in some cases and has exactly two accumulation points in other cases (Theorem 4.1). In case $\phi$ has infinite order, then the image set is even more complicated. We are able to prove that if $\phi(z)=\e^{(k)}(z)dz^2,$ for some positive integer $k$, where $\e^{(k)}(z)$ is defined inductively by   $\e^{(0)}(z)=1$ and $\e^{(j)}(z)=\e(\e^{(j-1)}(z))$, then  $ \overline{h(\C)}\cap\partial\H^2 =\cup_{j=0}^k\Cal A_j$ such that
   $A_j$  is countable and discrete for each  $0\le j\le k-1$,   $\Cal A_j$ consists of all isolated accumulation points of $\Cal A_{j-1}$ for $1\le j\le k$, and   $\Cal A_k$ consists of only one point (Theorem 4.2). 

In order to prove these results, we have to  study the regions where $|\phi|$ grows very fast and the regions where  $|\phi|$ decays or is bounded. In the regions where $|\phi|$ grows very fast, we refine the technique  in \cite{Hn, HTTW} which was introduced by Wolf and Minsky \cite{Wf, M}. In order to study the regions where $\phi$ is bounded, we need other tools. We will   use   the idea of the so-called maximal $\Phi$-radius of a holomorphic quadratic differential $\Phi$, see \S2 for definitions. Let $\Phi$ be  a holomorphic quadratic differential on $\D$. It was proved by Ani\'c, Markovi\'c and Mateljevi\'c \cite{A-M-M} that the norm of $\Phi$ with respect to the Poincar\'e metric is  uniformly bounded  if and only if the maximal $\Phi$-radius is uniformly bounded. On the other hand, it was proved in \cite{Wn} that  $h$ is a quasi-conformal harmonic diffeomorphism from $\H^2$ onto itself if and only if the norm of its Hopf differential is uniformally bounded. Hence we can conclude that $h$ is quasi-conformal if and only if the maximal $\Phi$-radius is uniformly bounded where $\Phi$ is the Hopf differential of $h$. In fact, it was proved in  \cite{A-M-M} that if $h$ is quasi-regular harmonic map on  $\H^2$ then the maximal $\Phi$-radius is uniformly bounded.  In this work, we will give a local version of these results. In Theorem 2.1, we will prove that if $h$ is an orientation preserving harmonic diffeomorphism from $\D$ or $\C$ into $\H^2$, under certain conditions, $h$ is quasi-conformal on  the domains where the maximal $\Phi$-radius is uniformly bounded, where $\Phi$ is the Hopf differential of $h$. In particular, we give another proof of  the result in \cite{Wn} mentioned above. Roughly  speaking,  in the case of harmonic diffeomorphisms  from $\C$ into $\H^2$, the domains where the $\Phi$-radius is uniformly bounded are the domains where $|\phi|$ is uniformly bounded and decays rapidly at infinity. 

In the process of proving Theorem 2.1, we need a refined version of the result in \cite{A-M-M} on the relation between the maximal $\Phi$-radius and the norm of the holomorphic quadratic differential  $\Phi$. In particular, we obtain a pointwise lower bound of the maximal $\Phi$-radius (Proposition 2.1). It turns out that the result also has applications to the problem of finding quasi-conformal harmonic diffeomorphism on $\H^2$ with prescribed quasi-symmetric function on the unit circle $\Bbb S^1$ which is identified as $\partial \H^2$. Let $\text {BQD}(\H^2)$ be the space of holomorphic quadratic differentials $\Phi$ on $\H^2$ such that
$$
|||\Phi|||=\sup_{z\in \H^2}||\Phi||(z)<\infty
$$
where $||\Phi||(z)$ is the norm of $\Phi$ at $z$ with respect to the Poincar\'e metric. In \cite{Wn}, a map $\b$ from $\bqd$ to the universal Teichm\"uller space $\t$ by sending $\Phi$ to the class of quasi-symmetric homeomorphism containing the boundary value of $h$, which is the quasi-conformal harmonic diffeomorphism on $\H^2$ with $\Phi$ as the Hopf differential. The map is injective \cite{L-T 3, L-W 2}  and  an open question is whether this map is surjective. This is in fact a conjecture of Schoen \cite{S}. There are partial results for this problem as well as similar problems in higher dimensions \cite{Ak, L-T 1--3,T-W 2-3,H-W, Y, S-T-W}. In our case, it is not hard to see that if one can prove that $\b$ is 'proper', namely, the inverse image of bounded set is bounded, then one can conclude that $\b$ is onto. Using the pointwise estimate of the maximal $\Phi$-radius and the main inequality of Riech and Strebel \cite{R-S}, we obtain  sufficient conditions for certain subspaces of $\bqd$ on which $\b$ is proper. For compact Riemann surfaces or Riemann surfaces of finite type, this kind of phenomena was   studied by Wolf \cite{Wf} and Markovi\'c-Mateljevi\'c \cite{M-M}. In \cite{M-M}, a generalized version of the inequality in \cite{R-S} was used. 

We organize the paper as follows. In \S1, we discuss some non-surjectivity results of harmonic maps from $\C$ into $\H^2$. In \S2, we study the relation between maximal $\Phi$-radius of the Hopf differential $\Phi$ of a harmonic map and quasi-conformality. In \S3 and \S4, we study the structures of images of harmonic maps from $\C$ into $\H^2$. In \S5, we use the result in \S2 to study quasi-conformal harmonic diffeomorphisms on $\H^2.$ In the appendix, we use Mathematica to produce figures of horizontal trajectories defined by different types of holomorphic quadratic differentials discussed in this work, so that one may get some feeling about the images of related harmonic maps.

Finally, the authors would like to thank the referee for pointing out a gap in the proof of theorem 1.1, which has been corrected accordingly.
\bigskip
\subheading{\S1 Results on non-surjectivity of harmonic diffeomorphisms}

In \cite{HTTW}, it was proved that a polynomial growth harmonic diffeomorphism from $\C$ into $\H^2$ is not surjective. In \cite{L-W 1}, the result was generalized to higher dimensions for polynomial growth harmonic maps between a more general class of manifolds. Not very many results are known if the map grows faster than polynomial. In this section, we will give  results on non-surjectivity of certain harmonic diffeomorphisms from $\C$ into $\H^2$ with fast growth rate. Note that   the growth rate of a harmonic diffeomorphism  from $\C$ into $\H^2$ 
can be expressed in terms of the growth rate of its Hopf differential, see \cite{T-Wn 1}. In particular, such a map is of polynomial growth if and only if its Hopf differential is of the form $Pdz^2$ with $P$ to be a polynomial.
  
\proclaim{Lemma 1.1} Let $\Omega$ be a domain in $\C$ which contains every disk  $\D( \i y, R(y))$ with center $\i y$ and radius $R(y)\ge 2\sqrt2(1+\epsilon)\log y $ for all $y\ge y_0>0$, where $\epsilon>0$ is a constant. Suppose $h$ is an orientation preserving  harmonic diffeomorphism from $\Omega$ into $\H^2$ with Hopf differential $\Phi=dz^2$. Then the length of the image of the half line $\Im z\ge y_0$, $\Re z=0$ under $h$ is bounded by a constant depending only on $\epsilon$ and $y_0$.
\endproclaim
\demo{Proof}  Let $\e(w)=||\partial h||$ be the norm of $\partial h$ and let $e$ be the energy density of $h$ with respect to the Euclidean metric in the domain, then the pull-back metric under $h$ is given by  
$$
h^*(ds_{\H^2}^2)=(e+2)dx^2+(e-2)dy^2=2\lf(\cosh(2w)+1\ri)dx^2+2\lf(\cosh(2w)-1\ri)dy^2.\tag1.1
$$
As in \cite{Wf, My} and page 63 in \cite{Hn} we can prove that there is $y_0>0$ such that if $y\ge y_0$
$$
0<w(\i y)\le C_1\e{\lf(-\frac{R(y)}{2\sqrt2}\ri)}\tag1.2
$$
 where $C_1$ is an  absolute constant.  Hence the length $\ell$ of the image of $\{\Im z\ge y_0,\, \Re z=0\}$ under $h$ satisfies:
$$
\split
\ell&=\int_{y_0}^\infty \lf[2\lf(\cosh(2w)-1\ri)\ri]^\frac12 dy\\
&\le C_3\int_{y_0}^\infty \e{\lf(-\frac{R(y)}{2\sqrt2}\ri)}dy\\
&\le C_4\int_{y_0}^\infty y^{-1-\epsilon}dy\\
&=C_5 \endsplit
$$
where $C_3--C_5$   are constants  depending only on $\epsilon$ and $y_0$, and we have used (1.1), (1.2) and the assumption that $R(y)\ge 2\sqrt2(1+\epsilon)\log y $ if $y\ge y_0$. The lemma then follows.
\enddemo
The following lemma  basically says that if $Q(z)=\frac12 z+o(1)$ as $\Re z\to\infty$, then the behavior of $\int \e(Q(z))dz$ is similar to that of $\int\e(\frac12z)dz$. 

\proclaim{Lemma 1.2}  Let $2\pi\ge A>0$ and let $Q(z)$ be an analytic function on the half strip 
$$\Cal S=\{z|\ \Re z>\alpha>0\ \text{and}\ \theta-A<\Im z<\theta+A\}$$ 
where $\theta$ is   constant. Suppose  $Q(z)=\frac12z+q(z)$ such that $|q(z)|\le g(\Re z)$ where $g(t)\ge0$ is a function defined on $\infty>t\ge\alpha$ which satsifies $\lim_{t\to\infty}g(t)=0$. Then for any $\frac14A>\delta>0$ there exists $\frak a>0$ depending only on $A$, $\delta$, $\alpha$ and the function $g$ such that if $z_0=x_0+\i \theta$ with $x_0>\frak a$ and if 
$$
\z(z)=\int_{z_0}^z\e\lf(Q(\xi)\ri)d\xi,
$$ 
then $\z$ maps $\Cal S_\delta$ injectively into $\z$-plane, and $\z(\Cal S_\delta)\supset \Cal R_{2\delta}\supset \z(\Cal S_{4\delta})$.
 Here
$$
\Cal S_{\delta}=\{z\in \Cal S|\ \Re z>x_0+\delta,\ \theta-A+\delta<\Im  z<\theta+A-\delta\},
$$
 $\Cal S_{4\delta}$ is defined similarly and
$$
\split
\Cal R_{2\delta}&=-2\exp(\frac{1}{2}z_0)+\lf\{\z|\ |\z|>2\e\lf(\frac12(x_0+2\delta)\ri),\right. \\ &\qquad\qquad \left. \frac12(\theta-A+2\delta)<\arg\z<\frac12(\theta+A-2\delta)\ri\}
.
\endsplit
$$
\endproclaim

\demo{Proof}
Since $\lim_{t\to\infty}g(t)=0$, for any $\epsilon>0$, there is $\frak a>\alpha$ depending only on $\alpha$ and $g$ such  that if $\Re z>\frak a$, then
$$
\lf|\e\lf(Q(z)\ri)-\e(\frac z2)\ri|\le \epsilon \e(\frac 12 \Re z).\tag1.3
$$
Let $x_0>\frak a$ and let $f(z)=2\lf(\e(\frac12 z)-\e(\frac12 z_0)\ri)$ with $z_0=x_0+i\theta$. Let $z_1=x_1+\i y_1$, $z_2=x_2+\i y_2$ in $\Cal S$ such that $x_1$, $x_2$ are larger than $x_0$. Suppose $x_1> x_2$, then by (1.3)
$$
\split
&\lf|\z(z_1)-\z(z_2)-f(z_1)+f(z_2)\ri|\\
\ &=|z_1-z_2|\lf|\int_0^1\lf(\e\lf(Q(tz_1+(1-t)z_2)\ri)-\e(\frac12(tz_1+(1-t)z_2)\ri)dt\ri|\\
&\le  \epsilon\lf(|x_1-x_2| +4\pi\ri) \e(\frac12 x_2)\int_0^1\e(\frac t2(x_1-x_2))dt\\
&\le \epsilon \e(\frac12 x_2)\lf[2\lf(\e(\frac12(x_1-x_2))-1\ri)+4\pi\e\lf(\frac12(x_1-x_2)\ri)\ri]\\
&\le  C_1\epsilon \e(\frac12 x_1)
\endsplit $$
where $C_1$ is an absolute  constant. Obviously, the inequality is still true if $x_1=x_2$. Hence, we have
$$
\lf|\z(z_1)-\z(z_2)-f(z_1)+f(z_2)\ri|\le C_1\epsilon\e(\frac12\max\{\Re z_1,\Re z_2\})\tag1.4
$$
provided $\Re z_1,\ \Re z_2>x_0$. 

On the other hand, for any $0<\delta_1<A$ and $z_1\neq z_2\in \Cal S$, with $\Re z_1\ge \Re z_2$ and $|\Im (z_1-z_2)|\le 2A-2\delta_1$, we have
$$
\split
|f(z_1)-f(z_2)|&=2\e(\frac12\Re z_1)|1-\e(\frac12(z_2-z_1))|\\
&\ge\tau \e(\frac12\Re z_1)
\endsplit$$
where $\tau>0$ depends only on $A-\delta_1$ and  the lower bound of $|z_1-z_2|$ where  we have used the fact that  $|\Im (z_1-z_2)|\le 2A-2\delta_1\le 4\pi-2\delta_1$. Hence for any $z_1\neq z_2\in \Cal S$,
$$
|f(z_1)-f(z_2)|\ge \tau\e(\frac12\max\{\Re z_1,\Re z_2\})\tag1.5
$$
where $\tau> 0$ depending only on the lower bound of $|z_1-z_2|$.
 
Let $0<\delta_1<A$, for any $a$ in $\Cal S_{\delta_1}\cap \{z|\ \Re z<\beta\}$ and $z$ the boundary of $\Cal S_{\frac12\delta_1}\cap\{z|\ \Re z<\beta+\frac12\delta_1\}$ where $\beta$ is a large number, we have
$$
\split
|f(z)-f(a)|&\ge \tau\e(\frac 12\Re z)\\
&\ge \frac{\tau}{C_1\epsilon}|\z(z)-f(z)|
\endsplit
$$
by (1.4) with $z_1=z$ and $z_2=z_0$, and (1.5), where $\tau>0$ is a constant depending only on $\delta_1$. Here we take $z_1=z$ and $z_2=z_0$ in (1.4). Choose $\epsilon$ small enough depending only on $A$ and $\delta_1$ such that $\frac{\tau}{C_1\epsilon}>1$, we have
$$
|f(z)-f(a)|>|\z(z)-f(z)|.
$$
Apply the Rouch\'e Theorem to the functions $\z-f(a)$, $f-f(a)$ on $\Cal S_{\frac12\delta_1}\cap\{z|\ \Re z<\beta+\frac12\delta_1\}$ and then let $\beta\to\infty $ we conclude that for any $a\in \Cal S_{\delta_1}$ there is one and only one $z\in S_{\frac12\delta_1}$ such that $\z(z)=f(a)$. 

On the other hand for such an $a$, we have
$$
\split
|\z(z)-\z(a)|&\ge |f(z)-f(a)|-|\z(z)-\z(a)-f(z)+f(a)|\\
&\ge \frac12|f(z)-f(a)|\\
&>|\z(z)-f(z)|\endsplit
$$
provided $\epsilon$ is chosen to be small enough (depending only on $A$ and $\delta_1$). Hence there is also exactly one   $z\in S_{\frac12\delta_1}$ such that $f(z)=\z(a)$. From these the lemma follows by considering the image of $f$.\enddemo

In the next theorem we will study the surjectivity of  those harmonic diffeomorphisms from $\C$ into $\H^2$ whose Hopf differentials are of finite order.  
\proclaim{Theorem 1.1} Let $h$ be an orientation preserving harmonic diffeomorphism from $\C$ into $\H^2$ with Hopf differential $\Phi=P\e(Q)dz^2$ such that
\roster
\item"{(i)}" $Q(z)=z^n+\sum_{j=1}^{n}a_jz^{n-j}$ is a polynomial of degree $n\ge 1$;
\item"{(ii)}" $P$ is entire with order $\rho<n$;
\item"{(iii)}"  there exists $\frac\pi n>\delta>0$ and $R_0>0$ such that 
$$
\Sigma\cap\{z|\ |z|>R_0 \ \text{and}\ -\delta<\arg z <\delta\}=\emptyset,
$$
where  $\Sigma$ is the set of all zeros of $P$.
\endroster
Then $h$ is not surjective. In particular, if $P$ is a polynomial then $h$ is not surjective.
\endproclaim
\demo{Proof} By the Hadamard factorization theorem, $P(z)=z^me^{a(z)}A(z)$, where $m$ is the order of $z=0$, $a(z)$ is a polynomial of degree less than $\rho<n$, and $A(z)$ is a canonical product of order less than or equal to $\rho$ formed by the zeros of $P$. So we can absorb $a(z)$ to the lower order terms of $Q(z)$ and assume $P(z)$ has the form $z^mA(z)$.

Let $\z_1=z^n$, which will map $-\delta <\arg z<\delta$ bijectively onto $-n\delta<\arg\z_1<n\delta$. In the region $$\Cal R_1=\{|\z_1|>R_0^n\}\cap\{-n\delta<\arg\z_1<n\delta\}.$$
$$
\Phi =n^{-2}\z_1^{-\frac{2(n-1)}n}\tilde P(\z_1)\e(\z_1)d\z_1^2
$$
where $\tilde P(\z_1)=P(z(\z_1))\e(\sum_{j=1}^na_j\z_1^{1-j/n})$. By (ii), without loss of generality we may assume that for $|\z_1|>R_0^n$, $|\tilde P(\z_1)|\le \e(|\z_1|^{\epsilon})$ for some $\epsilon>0$ which is small enough such that $\epsilon<1$. 
By (iii) and lemma 2.6.18 in \cite{B}, we have for any $\eta>0$, $\log|A(z)|>-|z|^{\rho+\eta}$ on $\{z|\ |z|>R_0 \ \text{and}\ -\delta<\arg z <\delta\}$ for a possibly larger $R_0$ and a smaller $\delta$. Therefore, we have
  $ \lf|\log|\tilde P(\z_1)|\ri| =O(|\z_1|^\epsilon)$ as $\z_1\to\infty$ for some $\epsilon<1$. From this, we have  $\lf|\nabla\log|\tilde P(\z_1)|\ri|=o(1)$ on $\tilde{\Cal R}_1={\Cal R}_1\cap \{-n\delta+\delta_1<\arg\z_1< n\delta-\delta_1\}$ for any $\delta_1>0$; and hence $|\frac{d}{d\z_1}\log\tilde P|=o(1)$ as $\z_1\to\infty$ and $\z_1\in \tilde{\Cal R}_1$. We conclude that
in $\tilde{\Cal R}_1$, 
$$
\Phi=\e(Q_1(\z_1))d\z_1^2=\phi d\z_1^2 
$$
where 
$$
Q_1(\z_1)=\z_1+Q_2(\z_1)
$$
with $-n\delta+\delta_1 <\arg\z_1<n\delta-\delta_1$ and $|\z_1|>R_0^n$. Here $Q_2(\z_1)=o(|\z_1|)$  and $\frac{d}{d\z_1}Q_2(\z_1)= o(1)$ as $\z_1\to\infty$. 

Let $\z_2=\z_1+Q_2(\z_1)$. This will map $\{-\pi <\arg\z_1<\pi \}\cap\{|\z_1|>R_1\}$ injectively onto its image for some $R_1>R_0^n$. Moreover,   there exists $R_2>0$ such that 
$\Cal R_2=\{-\pi+2\delta_1<\arg\z_2<\pi-2\delta_1\}\cap\{|\z_2|>R_2\}$ is in the image of the map $\z_1\mapsto \z_2$. In $\Cal R_2$, $\Phi$ can be written in the form
$$
\Phi=\lf[1+\frac{d  Q_2}{d\z_1}\ri]^{-2}\e(\z_2)d\z_2^2=\e(Q_3(\z_2))d\z_2^2
$$
where $Q_3(\z_2)=\z_2+o(1)$ as $|\z_2|\to\infty$. Let $\z(\z_2)=\int^{\z_2}\e(Q_3(\xi))d\xi$. By Lemma 1.2, we conclude that $\z(\z_2)$ will map a subdomain of $\Cal R_2$ bijectively onto the region
$$
\Cal R=\{\z|\ |\z|>R\ \text{and}\ \frac12\pi-\frac12\delta_2<\arg \z<\frac12\pi+\frac12\delta_2\}
$$
for some $R>0$ and $\delta_2>0$. On $\Cal R$, $\Phi=d\z^2$. The map $\z\mapsto\z_2\mapsto\z_1\mapsto z$ is injective when restricted on $\Cal R$. Hence $h(z(\z))$ is an orientation harmonic diffeomorphism from $\Cal R$ into $\H^2$. By Lemma 1.1, we conclude that the length of the image of the half line $\Im \z>a_0$, $\Re \z=0$ under $h$ is finite. Here $a_0$ is a large constant. By the definition of $\z$, 
$\Im \z\to\infty$ with $\Re \z=0$ implies that $z\to\infty$. Hence $h$ cannot be surjective. 
\enddemo

Please see the appendix for figures showing the behaviour of the horizontal trajectories for some typical examples of holomorphic quadratic differentials discussed here (figures 1-6). 
If we refine the method of proof in Theorem 1.1, we can generalize the result to some cases that the Hopf differentials grow very fast (see figure 7 in appendix). First we have the following:

\proclaim{Lemma 1.3}
Let $R_0>0$ and $\delta>0$ be constants and
let $h$ be an orientaion preserving harmonic diffeomorphism from $\Omega_{\delta}=\{ |z|> R_0,\, |\arg z|< \pi -\delta\}$ into $\H^2$ with Hopf differential of the form
$$
\Phi=\exp\left[ g_1+\exp\left[ g_2 +\cdots +\exp[g_k+Q]  \cdots \right]  \right] dz^2
$$
where $Q(z)=z^n+\sum_{j=1}^na_jz^{n-j}$ is a polynomial in $z$ and for each $j=1,\ldots,k$, $|g_j(z)|=O (\log |z|)$ as $|z|\to\infty$. Then there exists a path in $\Omega_{\delta}$ diverging to infinity such that its image under $h$ has finite length.
\endproclaim
\demo{Proof}
We will prove by induction on $k$. For $k=1$, $\Phi=\exp [g_1+Q]dz^2$ and we can apply the same proof as in Theorem 1.1 to conclude the existence of such path.

For $k\ge 2$, we consider the map $\zeta=f(z)=Q(z)+g_k(z)$ on a convex subdomain ${\Cal R}$ in $\Omega_\delta$ defined by
$$
\Cal R=\{z\in \Omega_\delta|\ \Re z\ge R_1,\  -\frac{\pi}{2(n-1)}+\epsilon<\arg z< \frac{\pi}{2(n-1)}-\epsilon\}
$$
 where $R_1>R_0$ and $\epsilon>0$ will be chosen later. It is clear that for $z=re^{i\theta}\in {\Cal R}$  with $r$ sufficiently large,  
$$
\Re f'(z)= nr^{n-1}\cos[(n-1)\theta]+o(r^{n-1})>0.
$$
Therefore, if we choose $R_1$ sufficiently large, $\zeta=f(z)$ maps ${\Cal R}$ one-one onto its image $f({\Cal R})$ as $\Cal R$ is convex (proposition 1.10 in \cite{P}). Since
$$
f(re^{\pm i(\frac{\pi}{2(n-1)} -\epsilon)})=
r^ne^{\pm i(n\frac{\pi}{2(n-1)} -n\epsilon)} +o(r^n),
$$
it is clear that $f({\Cal R})$ contains a subset of the form 
$$
\{ |\arg \zeta|<n\pi/2(n-1)-\epsilon_1,\, |\zeta|> R_2\}.
$$
If we choose $\epsilon<\frac{\pi}{2n(n-1)}$, then we can choose accordingly an $\epsilon_1$ such that $n\frac{\pi}{2(n-1)} -\epsilon_1 >\pi/2$. Hence, under this choice of $R_1$ and $\epsilon$,  $f({\Cal R})$ contains a half-plane $\{ \Re\zeta> R_2\}$ for some $R_2>0$.

On $f({\Cal R})$, in particular on the half-plane $\{ \Re\zeta> R_2\}$, $\Phi$ can be written as
$$
\Phi = \exp\left[\tilde{g}_1(\z) +\exp\left[ \tilde{g}_2(\zeta) +\cdots +\exp[\tilde{g}_{k-1}(\zeta)+\exp\zeta] \cdots \right]  \right] d\zeta^2 ,
$$
where $\tilde{g}_1(\zeta)=g_1(f^{-1}(\zeta))-2 \lf(\frac{d}{dz} \log f\ri)\lf (f^{-1}(\log\xi)\ri)$ and $\tilde{g}_j(\zeta)=g_{j}(f^{-1}(\zeta))$ for $j=2,\ldots, k-1$. 

Now, for any small $\delta_1>0$, lets consider the half strip
$$
{\Cal S}=\{ \Re\zeta >R_2,\, |\Im \zeta|<\pi-\delta_1\}.
$$
The exponential map $\xi=\exp\zeta$ maps ${\Cal S}$ one-one onto the domain
$$
\Omega_{\delta_1}=\{ |\xi| >R_3,\, |\arg \xi|<\pi-\delta_1\},
$$
where $R_3=e^{R_2}$. And on $\Omega_{\delta_1}$, $\Phi$ takes the form
$$
\Phi= \exp\left[ \tilde{\tilde{g}}_1(\xi)+\exp\left[ \tilde{\tilde{g}}_2(\xi) +\cdots +\exp[\tilde{\tilde{g}}_{k-1}(\xi) +\xi] \cdots \right]  \right] d\xi^2,
$$
where
$$
\tilde{\tilde{g}}_1(\xi)=g_1(f^{-1}(\log\xi))-2\lf(\frac{d}{dz} \log f\ri)\lf (f^{-1}(\log\xi)\ri)-2\log\xi
$$
and $\tilde{\tilde{g}}_j(\xi)=g_j(f^{-1}(\log\xi))$ for $j=2,\ldots,k-1$.

Since $|g_j(z)|=O(\log|z|)$ and $\zeta=f(z) \sim z^n$ as $|z|\to\infty$, we see that 
$$
|g_j(f^{-1}(\log\xi))|=O (\log\log|\xi|) \quad \text{for $j=1,\ldots, k-1$.}
$$
We also conclude from 
$\frac{df}{dz}(z)\sim z^{n-1}\sim \zeta^{(n-1)/n}$ as $|\zeta|\to\infty$ that 
$$
|\lf(\frac{d}{dz} \log f\ri)\lf (f^{-1}(\log\xi)\ri))|=O(\log\log|\xi|).
$$
All together we have, as $|\xi|\to\infty$,
$$
|\tilde{\tilde{g}}_j(\xi)|=\tilde{O}(\log|\xi|).
$$
Therefore, induction hypothesis implies that there exists a divergent path $\xi=\gamma(t)$ in $\Omega_{\delta_1}$ such that its image under the harmonic maps $h\circ f^{-1}\circ \log$ has finite length. That is, there exists a path $f^{-1}(\log\gamma(t))$ in $\Omega_{\delta}$ such that its image under $h$ has finite length.
\enddemo

\proclaim{Theorem 1.2} Let $h$ be a harmonic diffeomorphism from $\C$ into $\H^2$ with Hopf differential of the form
$$
\Phi=P_1\e\left[P_2\e\left[\cdots \e\left[P_k\e(Q)\right]\cdots\right]\right]dz^2, 
$$
where $Q(z)=a_nz^n+\cdots$ and $P_j$, $j=1,\ldots,k$ are polynomials and $k\ge 1$.
 Then $h$ is not surjective.
\endproclaim
\demo{Proof} By making a change of parameter of the form $z\to r_0\e(\i \theta_0)z$ for some constants $r_0$ and $\theta_0$, we may assume that $Q=z^n+\sum_{j=1}^na_jz^{n-j}$. As $P_j$, $j=1,\ldots,k$, are polynomials, there exists $R_0>0$ such that there is no zeros of any $P_j$ in the set $\{|z|>R_0\}$. Then for any $\delta>0$, one can define $g_j=\log P_j$ on the set $\Omega_{\delta}=\{ |z|>R_0,\, |\arg z|<\pi-\delta\}$ and the Hopf differential $\Phi$ can be written in the form required in Lemma 1.3. Therefore, there exists a path diverging to infinity in $\Omega_{\delta}$ such that its image under $h$ has finite length. Hence, $h$ is not surjective.
\enddemo

Let us finish the dicussion of this section by given a different type of condition for the nonsurjectivity.
Recall that a complex number $a$ in the extended complex plane is said to be an {\it asymptotic value} of an entire function $f(z)$ if there is a path $z(t)$, $0\le t<1$ such that $\lim_{t\to 1}z(t)=\infty$ and $\lim_{t\to 1}f(z(t))=a$. If $a$ is a finite number, then it is called a finite asymptotic value. 

\proclaim{Theorem 1.3} Let $f$ be an entire function.  Suppose there exist $\delta>0$ and $R>0$ such that
\roster
\item"{(i)}" $f$ has no finite asymptotic value in the domain
$$
\Cal R=\left\{\zeta|\ \frac\pi2-\delta<\arg\z<\frac\pi2+\delta,\ \text{and}\ |\z|>R\right\};\ \text{and}
$$

\item"{(ii)}" $f'(z)\neq0$ for all $z$ in $f^{-1}(\Cal R)$.
\endroster
Suppose $h$ is an orientation harmonic diffeomorphism from $\C$ into $\H^2$ with Hopf differential $(f')^2dz^2$, then $h$ is not surjective.
\endproclaim
\demo{Proof} Let $\Omega$ be a component of $f^{-1}(\Cal R)$. By (ii), $f$ is a local diffeomorphism on $\Omega$. By (i), we can conclude that every path in $\Cal R$ begins at $\z_0$ can be lifted to a path in $\Omega$ which begins at a point $z_0$ with $f(z_0)=\z_0$. Since $\Cal R$ is simply connected, $f$ maps $\Omega$ bijectively to $\Cal R$. Hence $h\circ f^{-1}(\z)$ is a harmonic diffeomorphism from $\Cal R$ into $\H^2$ with Hopf differential $d\z^2$. Moreover, $f^{-1}(\z)\to\infty$ if $\z\in \Cal R$ and $\z\to\infty$. The result follows from Lemma 1.1.

\enddemo
\bigskip
\subheading{\S2 Maximal $\Phi$-radius  and quasi-conformal harmonic maps}

Let us recall the definition of maximal $\Phi$-radius of a holomorphic quadratic differential $\Phi$ on a domain in $\C$. Let $\Omega$ be a domain in $\C$ and let $\Phi=\phi dz^2$ be a holomorphic quadratic differential on $\Omega$. Let $z_0\in \Omega$ such that $\phi(z_0)\neq0$. Choose a branch of $\sqrt{\phi}$ near $z_0$, and let 
$$
w=f(z)=\int_{z_0}^z\sqrt{\phi(\z)}d\z.
$$
Let $B(R)=\{w|\ |w|<R\}$  be the  maximal disk in the $w$-plane such that $f^{-1}$ is a conformal diffeomorphism from $B(R)$ into $\Omega$. Then $R_{z_0,\Omega}=R$ is called the {\it maximal $\Phi$-radius} of $\Phi$ at $z_0$ with respect to $\Omega$ and $V_{z_0,\Omega}=f^{-1}(B(R))$ is called the {\it maximal $\Phi$ disk} around $z_0$ with respect to $\Omega$. We will drop the subscript $\Omega$ if this will not cause any confusion. Moreover, by convention if $\phi(z_0)=0$ we define $R_{z_0}=0$.  In \cite{Hn}, it was proved that if $h$ is an orientation preserving harmonic diffeomorphism from a domain $\Omega$ in $\C$ into $\H^2$ with Hopf differential $\Phi$, and if $z_n$ is a sequence in $\Omega$ with $R_{z_n}\to\infty$, then the modulus the complex dilatation 
of $h$ at $z_n$ will tend to 1. Conversely, one would like to know whether $h$ would be quasi-conformal on a set with bounded maximal $\Phi$-radius. In this section, we will prove that this is the case under certain assumptions. The result will be useful to study images of harmonic diffeomorphisms from $\C$ into $\H^2$.

Let $\Omega$ be a hyperbolic domain in $\C$, i.e.   its universal cover is conformal to the unit disk. Let $\rho^2|dz|^2$ be the hyperbolic metric on $\Omega$, i.e.  the complete metric  with constant Gaussian curvature $-1$.  Then it is known that \cite{Ah} for any $z\in\Omega$
$$
\rho(z)\le \frac{2}{d(z,\partial \Omega)},
$$
where $d(z,\partial \Omega)$ is the Euclidean distance from $z$ to $\partial \Omega$. If in addition, we have
$$
\rho(z)\ge \frac{C}{d(z,\partial \Omega)}
$$
for some positive constant $C$ for all $z\in \Omega$, then we say that $\Omega$ is   {\it strongly hyperbolic}. Please note that our definition is slightly different from that in \cite{A-M-M}. It is shown in Theorem 5 of \cite{A-M-M}  that if $\Omega$ is bounded hyperbolic and the diameters of the boundary  components are uniformly bounded from below by a positive constant then $\Omega$ is strongly hyperbolic.  Moreover, being strongly hyperbolic is conformally invariant:

\proclaim{Lemma 2.1} Let $\Omega_1$ and $\Omega_2$ be conformally equivalent domains. Suppose $\Omega_1$ is strongly hyperbolic, so is $\Omega_2$. 
\endproclaim
\demo{Proof} Obviously $\Omega_2$ is hyperbolic. Let $w=f(z)$ be a conformal diffeomorphism from $\Omega_1$ onto $\Omega_2$.  Let $\rho_2^2|dw|^2$ be the hyperbolic metric on $\Omega_2$ and let $\rho_1^2|dz|^2$ be the hyperbolic metric on $\Omega_1$. Let $d_1(z)=\dist(z,\partial\Omega_1)$ and $d_2(w)=\dist(w,\partial\Omega_2)$, where both distances are   Euclidean distances. Then by well-known fact \cite{V, p. 147}, we have
$$
\frac14d_1(z)|f'(z)|\le  d_2(f(z)).
$$
 Hence
$$
\split
\rho_2(f(z))&=|f'(z)|^{-1}\rho_1(z)\\
&\ge \frac{C}{d_1(z)|f'(z)|}\\
&\ge \frac{C}4d_2(f(z))
\endsplit
$$
for some positive contant $C$, where we have used the fact that $\Omega_1$ is strongly hyperbolic. Hence $\Omega_2$ is strongly hyperbolic.
 \enddemo

 Let $\Omega$ be a strongly hyperbolic domain and let $\Phi=\phi dz^2$ be a holomorphic quadratic differential on $\Omega$ with hyperbolic metric $\rho^2|dz|^2$. For $z\in \Omega$, let $||\Phi(z)||=\rho^{-2}(z)|\phi|(z)$ be the norm of $\Phi$ at $z$ and let $|||\Phi|||=\sup_{z\in \Omega}||\Phi||(z)$. The following is proved in \cite{A-M-M} (the equation (2) and the lemma 1.2)
\proclaim{Theorem}(Ani\'c-Markovi\'c-Mateljevi\'c) With the above notations with $\Omega$ being the unit disk $\D$, there exists an absolute constant $C>0$ such that for any holomorphic quadratic differential $\Phi$ on $\D$ we have
$$
||\Phi||(z)\ge C^{-1}R^2_z\tag2.1
$$
for all $z\in \D$, 
and
$$
|||\Phi|||\le C R^2_\infty\tag2.2
$$
where $R_\infty=\sup_{z\in\D}R_z$. 
\endproclaim

They actually proved that (2.1) is true for any hyperbolic domain in $\C$. We will  obtain a pointwise estimate  for strongly hyperbolic domain which implies (2.2). The estimate will be useful  in applications. 

\proclaim{Proposition 2.1}Let $\Omega$ be a hyperbolic domain and  let $\Phi=\phi dz^2$ be a holomorphic quadratic differential defined on $\Omega$. Then there exists an absolute positive constant $C$ such that for $z\in \Omega$
$$
 R_z\le C||\Phi||^{\frac12}(z).  \tag2.3
$$
If in addition $\Omega$ is strongly hyperbolic then there is a positive constant $C'$ depending only on $\Omega$ such that for $z\in \Omega$ with $\Phi(z)\neq0$
$$
R_z\ge C' \frac{||\Phi||^2(z)}{|||\Phi|||^\frac32}\tag2.4
$$
\endproclaim

We need  the following lemmas.
\proclaim{Lemma 2.2} Denote $B(r)$ to be the set of complex numbers with modulus less than $r$. Let $f: B(r)\to \C$ be an analytic function, such that
$f(0)=0$ and $f'(0)\neq0$. Suppose $|f(z)|\le M$ for all $z$, then (i) $f$ is
one to one on $B(r_1)$, where $r_1=\frac{r^2|f'(0)|}{8M}$; (ii) $f(B(r_1))\supset B(\frac{r^2|f'(0)|^2}{32M})$.
\endproclaim
\demo{Proof} Let us first assume that $r=1$, and $f'(0)=1$. Then 
$$
f(z)=z+\sum_{n=2}^{\infty}a_nz^n.
$$
By Cauchy theorem, we have $|a_n|\le M$, and $1\le M$. Suppose $z_1\neq z_2$ are
in
$B(\frac1{8M})$, and $r=\max\{|z_1|,\ |z_2|\}$, then
$$
\split
|f(z_1)-f(z_2)|&=\lf|(z_1-z_2)+\sum_{n=2}^{\infty}a_n(z_1^n-z_2^n)\ri|\\
&\ge|z_1-z_2|\lf|1-M\sum_{n=2}^{\infty}nr^{n-1}\ri|\\
&=|z_1-z_2|\lf|1-Mr\frac{2-r}{(1-r)^2}\ri|\\
&\ge |z_1-z_2| (1-8Mr)\\
&>0
\endsplit
$$
where we have used the facts that $M\ge1$, and $r<\frac1{8M}<\frac12$. Hence $f$ is
one to one on $B(\frac1{8M})$. Using  the fact that the Koebe's constant is
$\frac14$ \cite{V, p. 149}, we have
$f\lf(B(\frac1{8M})\ri)$ contains
$B(\frac1{32M})$. In general, if $f$ is defined on $B(r)$ with $f(0)=0$ and with
$b=f'(0)\neq 0$. Define $\tilde f(\zeta)=\frac{f(r\zeta)}{rb}$ for $\zeta\in
\D$. Then $\tilde f(0)=0$, and $\tilde f'(0)=1$. Let $M_1=\frac{M}{r|b|}$, then
$M_1\ge |\tilde f(\zeta)|$ for all $\zeta$. Hence $\tilde f$ is one to one on
$B(\frac{1}{8M_1})$ and $\tilde f\lf(B(\frac1{8M_1})\ri)$ contains
$B(\frac1{32M_1})$. Hence $f$ is one to one on
$B(\frac{r}{8M_1})=B(\frac{r^2|f'(0)|} {8M})=B(r_1)$, and
$f\lf(B(r_1)\ri)$ contains $B(\frac{r|b|}{32M_1})=B(\frac{r^2|f'(0)|^2}{32M})$.
\enddemo  

The following lemma is proved in \cite{A-M-M, see Lemma 1.2}.
\proclaim{Lemma 2.3} Let $f:B(r)\to \C$ be analytic, and $M\ge |f(z)|$ for all $z$.
Suppose $f(0)\neq0$, then $f(z)\neq0$ for all   $z\in B(\frac{r|f(0)|}{2M}).$
\endproclaim

\demo{Proof of Proposition 2.1} (2.3) was  proved in the Lemma 2.3 of \cite{A-M-M}. In order to prove (2.4), let $z_0\in \Omega$, with $ \phi(z_0) \neq0$. $\phi$
is analytic on $B(z_0,r_0)$ where $r_0=\frac 12 d(z_0,\partial \Omega)$, where $d$ is the Euclidean distance. By Lemma 2.3, $\phi$ is
never zero on $B(z_0, r)$, where 
$$r=\frac{r_0|\phi(z_0) |}{2M_0}
$$
and $M_0=\sup_{B(z_0,r_0)}|\phi|$. Hence we can take a branch of square root of $\phi$ in
$B(z_0,r)$. Let $f(z)=\int_{z_0}^z\sqrt \phi(\z)d\z$, for $z\in B(z_0,r)$, then $f$
is analytic, $f(z_0)=0$ and $|f'(z_0)|=|\phi(z_0)|^\frac12\neq0$. By Lemma 2.2, $f$
is one to one on $B(r_1)$ where $r_1=\frac{r^2|f'(z_0)|}{8M_1}$, where
$M_1=\sup_{B(z_0,r)}|f|$. Moreover, $f(B(r_1))$ contains the disk
$B(R)=B(\frac{r^2 |f'(z_0)|^2}{M_1})$. Now $M_1\le rM_0^\frac12 $implies that 
$$
R_{z_0}\ge R=\frac{r^2|\phi(z_0)|}{32rM_0^\frac12}=\frac{r_0|\phi(z_0)|^2}{64M_0^\frac32}.
\tag2.5
$$
This will imply the proposition because $\Omega$ is strongly hyperbolic.
\enddemo
From the proof of the proposition, we have the following corollary which will be used in \S3 and \S4.
\proclaim{Corollary 2.1} Suppose $\phi$ is analytic on $B_{z_0}( R)$ such that $\alpha|\phi|(z_0)\ge |\phi|(z) $ for some constant $\alpha>0$  for all $z\in B_{z_0}(R)$. Let $\Phi=\phi dz^2$. Then the maximal $\Phi$-radius of $z_0$ with respect to $B_{z_0}( R)$ is bounded below by $\frac{R|\phi|^\frac12(z_0)}{64\alpha^\frac32}$.
\endproclaim
\demo{Proof} This is a direct consequence of (2.5).
\enddemo

\proclaim{Lemma 2.4} Let $\Omega$ be a simply connected domain in $\C$ and $\Phi=\phi dz^2$ be a holomorphic quadratic differential on $\Omega$. Let $z_0\in \Omega$ such that $\phi(z_0)\neq0$ and let $R$ be the maximal $\Phi$-radius of $z_0$ with maximal $\Phi$-disk $V$. Suppose $R<\infty$ and suppose
$$
w=\psi(z)=\int_{z_0}^z\sqrt{\phi}d\z 
$$
$z\in V$. Then for any $0<\delta<R$, there exists a point $z\in V$ with $|\psi(z)|=\delta$ such that the $\Phi$-radius of $z$ is exactly $R-\delta$.
\endproclaim
\demo{Proof} By the definitions of $V$ and $R$,  $\psi^{-1}:\D_R\to V$ is a bijective conformal diffeomorphism, where
$$
\D_R=\{w|\ |w|<R\}.
$$
   It is easy to see that if $z\in V$ with $|\psi(z)|=\delta$ then the  $\Phi$-radius of $z$ is at least $R-\delta$. Suppose the lemma is not true, then $R_z>R-\delta$ for any $z\in\psi^{-1}\lf(\{w|\ |w|=\delta\}\ri)$. Because $R_z$ is continuous, there is $\epsilon>0$ such that  $R_z\ge R-\delta+\epsilon$ for all $z$ with $|\psi(z)|=\delta$. Hence $\psi^{-1}$ can be extended to an analytic function from $\D_{R+\epsilon}$ to $\Omega$ such that it is a local diffeomorphism. In particular, $\phi$ is not zero in   $\psi^{-1}(\D_{R+\epsilon})$.  By the definition of $R$, there exist two sequences  $w_n$ and $\tilde w_n$ such that for each $n$ both $w_n$ and $\tilde w_n$ are in $\D_{R+\epsilon/n}$, $w_n\neq\tilde w_n$ but $\psi^{-1}(w_n)=\psi^{-1}(\tilde w_n)$. Without loss of generality, we may assume that $w_n\to a$ and $\tilde w_n\to b$, and 
$$
\lim_{n\to\infty}\psi^{-1}(w_n)=\lim_{n\to \infty}\psi^{-1}(\tilde w_n)=c.
$$
Since $\psi^{-1}$ is a local diffeomorphism, $a\neq b$. Note that $a$ and $b$ are in $\overline\D_R$. Let $\gamma$ be the straight line joining $a$ and $b$ and let $\Gamma=\psi^{-1}(\gamma)$. Then $\Gamma$ is a smooth simple closed curve in $\Omega$ because $\psi^{-1}$ is one to one on $\D_R$. Let $\theta$ be the interior angle at $c$.  Apply the Gauss-Bonnet Theorem for the metric $(|\phi|+\eta )|dz|^2$ on $\Omega$ with $\eta>0$, we have
$$
-\frac12\int_{\Omega_1}\Delta\log(|\phi|+\eta)+\int_{\Gamma}\kappa_\eta=\pi+\theta 
$$
where $\Omega_1$ is the interior of $\Gamma$, $\kappa_\eta$ is the geodesic curvature of $\Gamma$ with respect to the metric $(|\phi|+\eta )|dz|^2$ and $\Delta$ is the Euclidean Laplacian. Here we have used the fact that $\Omega$ is simply connected. Let $c_1,\dots,c_\ell$ be the zeros of $\phi$ inside  $\Omega_1$  with multiplicities $k_1,\dots,k_\ell$ with $k_j\ge0$. Let $r>0$ be small enough so that  $1\le j\le \ell$ the disks $D_j$   of radius $r$ and centers at $c_j$ are disjoint and are  inside of $\Omega_1$. Then we have
$$
\split
\pi+\theta&=-\frac12\sum_{j=1}^\ell \int_{D_j}\Delta \log(|\phi|+\eta)- \frac12 \int_{\Omega_1\setminus\cup_{j=1}^\ell D_j}\Delta \log(|\phi|+\eta)+\int_{\Gamma}\kappa_\eta\\
&=-\frac12\sum_{j=1}^\ell \int_{\partial  D_j}\frac{\partial}{\partial r}\log(|\phi|+\eta)- \frac12 \int_{\Omega_1\setminus\cup_{j=1}^\ell D_j}\Delta \log(|\phi|+\eta)+\int_{\Gamma}\kappa_\eta\\
&\to -\frac12\sum_{j=1}^\ell \int_{\partial D_j}\frac{\partial}{\partial r}\log |\phi|+\int_{\Gamma}\kappa
\endsplit
$$
as $\eta\to0$, where $\kappa$ is the geodesic curvature with respect to the metric $|\phi||dz|^2$ and we have used the fact that $\log|\phi|$ is harmonic. Since $|\phi||dz|^2$ is the pull-back metric under $\psi^{-1}$ of the flat metric in the $\D_R$, we have $\kappa=0$ on $\Gamma$. Let $r\to0$, we conclude that
$$
\pi+\theta=-\sum_{j=1}^\ell k_j\pi.
$$
Since $\theta\ge0$ and $k_j\ge0$ for all $j$, this is impossible.
\enddemo

\proclaim{Theorem 2.1} Let $\Omega$ be a strongly hyperbolic domain in $\C$   with hyperbolic metric $e^{2v}|dz|^2$ and $\Omega_1\subset \Omega$. Let $h$ be an orientation preserving harmonic diffeomorphism from $\C$ into $\H^2$ and let $w=\log|\partial h|$, where $|\partial h|$ is the norm of $\partial h$ with respect to the Euclidean metric on $\Omega$ and hyperbolic metric on $\H^2$.
 Let $\Phi=\phi dz^2$ be the Hopf differential of $h$ and let $R_z$ be the maximal $\Phi$-radius of $z$ with respect to $\C$. Suppose
$\sup_{z\in \Omega}R_z=R<\infty$ and $w\ge v-C$ on $\Omega_1$ for some constant $C$ and $\inf_{z\in \C\setminus \Omega_1}R_z>0$. Then $h$ 
is quasi-conformal on $\Omega_1$.
\endproclaim
\demo{Proof}  Suppose that $h$ is not quasi-conformal on $\Omega_1$. Then there exists  $z_n\in \Omega_1$ such that $ \phi(z_n)e^{-2w(z_n)}\to 1$ as $n\to\infty$. Since $R_{z_n}\le R$,  we may assume that $\lim_{n\to\infty} R_{z_n}=R_0$. Suppose $R_0>0$. Let $V_{z_n}$ be the maximal $\Phi$-disk with image $\D_{R_{z_n}}$ and let $\z=\int_{z_n}^z\sqrt\phi dz=\psi_n(z)$. Let $\tilde w_n(\z)=(w -\frac12\log|\phi|)(\psi_n^{-1}(\z))>0$ which is considered as a function on $\D_{R_{z_n}}$. Then $\tilde w_n>0$ and
$$
\Delta_\z\tilde w_n=e^{2\tilde w_n}-e^{-2\tilde w_n}.
$$
Then $\tilde w_n$ are locally uniformly bounded by proposition 1.5 in \cite{T-W 1}. Passing to a subsequence if necessary, $\tilde w_n$ converges uniformly on compact subsets of $\D_{R_0}$. Since $\tilde w_n(0)\to 0$, by mean value inequality, $\tilde w_n\to 0$ uniformly on compact sets of $\D_{R_0}$. By Lemma 2.4,  for each $n$ there exists $\z_n$ with $|\z_n|=\frac12R_{z_n}$ such that the $\Phi$-radius of $z_n'=\psi_n^{-1}(\z_n)$ is $\frac12R_{z_n}$.  Moreover, we still have $\phi(z'_n)e^{-2w(z'_n)}\to 1$. Continue in this way  and by a diagonal process, if $h$ is not quasi-conformal on $\Omega_1$, then we can find $z_n\in \C$ such that 
$$
\lim_{n\to\infty}|\phi(z_n)|e^{-2w(z_n)}=1\tag2.6
$$
and
$$
\lim_{n\to\infty}R_{z_n}=0.\tag2.7
$$
Since $\inf_{z\in\C\setminus\Omega_1}R_z>0$, we may assume that $z_n\in\Omega_1$ for all $n$. By Proposition 2.1, we have
$$
\split
R_{z_n}&\ge C_3\frac{||\Phi||^2(z_n)}{|||\Phi|||^\frac32}\\
&\ge \frac{C_4}{R^3}\lf[|\phi(z_n)|e^{-2v(z_n)}\ri]^2\\
&\ge \frac{C_5}{R^3}\lf[|\phi(z_n)|e^{-2w(z_n)}\ri]^2,
\endsplit\tag2.8
$$
where the norm of $\Phi$ is taken with respect to the metric $e^{2v}|dz|^2$ on $\Omega$. Here we have used (2.3), (2.4),  the fact that the $\Phi$-radius with respect to $\Omega_1$ or $\Omega$ is no greater than the $\Phi$-radius with respect to $\C$,  the assumption that $R_z$ are uniformly bounded by $R$ on $\Omega$ and that $w\ge v-C$ on $\Omega_1$. Let $n\to\infty$ in (2.8), we have a contradiction because of (2.6) and (2.7). This completes the proof of the theorem.
 \enddemo

\proclaim{Remark 2.1} In the theorem, we may replace $\C$ by the unit disk. Moreover, suppose $\Omega$ is a subset of $\C$ (respectively $\H^2$) and $h$ is  an orientation preserving harmonic diffeomorphism from $\C$ (respectively $\H^2$) into $\H^2$ such that $|\partial h|^2|dz|^2$ is complete in $\C$ (respectively $\H^2$), where the norm   is taken with respect to the Euclidean metric in the domain. 
Then the assumption that $w\ge v-C$ on $\Omega_1$ in the theorem can be replaced by $w\ge v-C$ on $\partial \Omega_1$ by the comparison principle in \cite{Wn}.\endproclaim

By (2.3), which was proved in \cite{A-M-M}, and the above remark, we obtain a new proof of the following result in \cite{Wn} as a corollary of Theorem 2.1.

\proclaim{Corollary 2.2} Let $h$ be an orientation preserving harmonic diffeomorphism on $\H^2$ with Hopf differential $\Phi$ such that $|\partial h|^2|dz|^2$ is complete. Suppose $|||\Phi|||<\infty$, then $h$ is quasi-conformal. Here the norm of $\partial h$ is taken with respect to the Euclidean metric in the domain and the norm of $\Phi$ is taken with respect to the Poincar\'e metric while $|||\Phi|||$ is taken with respect to the Poincar\'e metric.
\endproclaim

\bigskip
\subheading{\S3 Image of harmonic diffeomorphism with Hopf differential $P\e(Q)dz^2$}

Let $h$ be an orientation preserving harmonic diffeomorphism from $\C$ into $\H^2$ with Hopf differential $\Phi=\phi dz^2=P\e(Q)dz^2$ where $P$ and $Q$ are polynomials. By the result of \S1, we know that $h$ is not surjective. Assume that $|\partial h|^2|dz|^2$ is complete on $\C$. In \cite{HTTW}, it was proved that if $Q$ is a constant, that is, if $\phi$ is a polynomial of degree $m$, then the closure of the image of $h$ in $\overline{\H^2}$ is the convex hull of an ideal polygon with $m+2$ vertices in $\H^2$. The result is generalized from $\C$ to surfaces with finite total curvature and in higher dimensions in \cite{L-W 1, 2}. The assumption that $\phi$ is a polynomial is equivalent to the fact that $h$ is of polynomial growth. Let $A$ be the intersection of the closure of the image of $h$ with the geometric boundary of $\partial \H^2$. If $Q$ is not constant, then $A$ will no longer be a finite set. In this section, we will prove that in this case, $A$ is a countable set with exactly $n$ distinct accumulation points where $n$ is the degree of $Q$. In fact, we will prove that the result is true for a larger class of harmonic diffeomorphisms. 

First we need a lemma. For $\alpha\ge0$, let 
$$
\C_{\alpha}=\{z|\  \Re z> \alpha, \ -\infty<\Im z<\infty\}.$$ 
Let $h$ be an orientation preserving harmonic diffeomorphism from $\C_\alpha$ into $\H^2$ with Hopf differential  $\Phi=\e(Q)dz^2$ with
$Q(z)= z+q(z)$ such that
$|q(z)|\le g(\Re z)$ for some nonnegative function $g$ with $\lim_{t\to\infty}g(t)=0$. 

\proclaim{Lemma 3.1} With the above notations and assumptions, we have the following:
\roster
\item"{(i)}" There exist distinct points $p_k\in \partial\H^2$ with $k=0,\pm1,\pm2,\dots$ such that for any $\pi>\delta>0$ and for any sequence $z_n\in\C_\alpha$ with $(2k-1)\pi+\delta<\Im z_n<(2k+1)\pi-\delta$ and $z_n\to\infty$, then $h(z_n)\to p_k$ as $n\to\infty$. Moreover, the $p_k$'s are monotone in $\Bbb S^1$.

\item"{(ii)}" For all $\delta>0$ and $\epsilon>0$, there is $a>0$ such that  for any integer $k$, if $2k\pi+\delta<\Im z<2(k+1)\pi-\delta$, then $d_{\H^2}(h(z),\gamma_k))\le \epsilon$ for all $z\in\C_{\alpha}$ with $\Re z>a$, where $\gamma_k$ is the geodesic joining $p_k$ and $p_{k+1}$.

\item"{(iii)}" There is $b>0$ such that if $z_n\in \C_\alpha$   $\Re z_n\ge b$ and $\Im z_n\to+\infty$ (respectively $\Im z_n\to-\infty$), then $\lim_{n\to\infty}h(z_n)=p_+$ (respectively $\lim_{n\to\infty}h(z_n)=p_-$), where
$p_+=\lim_{k\to\infty}p_k$ and $p_-=\lim_{k\to-\infty}p_k$.
\endroster
\endproclaim
\demo{Proof} For simplicity, assume $\alpha=0$. To prove the existence of those $p_k\in\partial\H^2$ in (i), we apply  Lemma 1.2 with $\theta=0$ and $A=2\pi$ to obtain that for any $\pi>\delta>0$, there exists $x_0>0$ such that if $z_0=x_0$ and $\z(z)=\int_{z_0}^z\e(\frac12Q(\xi))d\xi+\e(\frac12 x_0)$, then 
$\z$ is injective on $\Cal S_{\frac14\delta}$, and $\z(\Cal S_{\frac14\delta})\supset \Cal R_{\frac12\delta}\supset \z(\Cal S_\delta)$ where $\Cal S_\delta$ and $\Cal R_\delta$ etc. are defined as in Lemma 1.2 (with $A=2\pi$ and $\theta=0$). Then $h(z(\z))$ is an orientation preserving harmonic diffeomorphism from $\z(\Cal S_{\frac14\delta})$ into $\H^2$ with Hopf differential $\Phi=d\z^2$.  Note that the maximal $\Phi$-radius of any point $\z=u+\i v$ in $\Cal R_{\frac12\delta}$ is at least $u-C_1$ for some constant $C_1$ depending only on $\delta$. As in \cite{HTTW, p.109}, we can prove that  the image of any horizontal half line $ \z(t)=t+\i v_0 $ with $t $ being larger than some constant in $\Cal R_{\frac12\delta}$ under $h$ is asymptotically a geodesic near infinity and   tends to a point  in $\partial \H^2$ as $t\to\infty$. By the proof of Lemma 1.1, we can conclude that the image of any vertical line $u=$constant in $\Cal R_{\frac12\delta}$ under $h$ has uniformly bounded length. Hence if $\z_n\in  \Cal R_{\frac12\delta}$, $\z_n\to\infty$ then $h(z(\z_n))\to p_0$ for some $p_0\in \partial \H^2$. Since $\z(\Cal S_\delta)\subset  \Cal R_{\frac12\delta}$ and $z_n\to\infty$ implies that $\z(z_n)\to\infty$ for $z_n\in\Cal S_\delta$, we have
$$
\lim_{n\to\infty}h(z_n)=p_0.
$$
Similarly, one can prove that for any integer $k$ there exists $p_k\in \partial \H^2$ such that for   any $\delta>0$  and $z_n$ with   $ (2k-1)\pi+\delta<\Im z_n<(2k+1)\pi-\delta$ such that $z_n\to\infty$, then
$$
\lim_{n\to\infty}h(z_n)=p_k.
$$

To prove the remaining of (i) and (ii), we use  Lemma 1.2 again to conclude that for all $\delta>0$ small enough, there exist  $a_j>0$, $b_j>0$, $j=1,\ 2$ such that for any integer $k$, there is a analytic function $\z=\z^{(k)}(z)$ which maps 
$$
\Cal S_1=\{z|\ \Re z>a_1, \ (2k-1)\pi+\frac12\delta <\Im z<(2k+3)\pi-\frac12\delta\}
$$ 
and
$$
\Cal S_2=\{z|\ \Re z>a_2, \ 2k\pi- \delta <\Im z<2(k+1)\pi+ \delta\}
$$
injectively into $\z$-plane. Moreover, $\z(\Cal S_1)\supset \Cal R_1$, $\z(\Cal S_2)\subset\Cal R_2$, for $j=1, \ 2$. Here
$$
\Cal R_1=\{\z|\ |\z|>b_1,\ \frac12(2k-1)\pi+\delta<\arg\z<\frac12(2k+3)\pi-\delta\},
$$
$$
\Cal R_2=\lf\{\z|\ |\z|>b_2, \ \frac12\lf(2k\pi-\frac32\delta \ri) <\arg\z<\frac12\lf(2(k+1)\pi+\frac32\delta\ri)\ri\}
$$
for $j=1,\ 2$. Moreover, the Hopf differential $\Phi$ of $h$ in the $\z$ coordinates is of the form $d\z^2$. We will write $h(\z)$ instead of $h(z(\z))$ if   no confusion will arise. Let us consider the case when $k$ is even. The case that $k$ is odd is similar. By the previous result, we know that if  $\z_n\in \Cal R_1$ with $\Re \z_n\to\infty$ along a half line $\Im \z$=constant, then $h(\z_n)\to p_k$, and if $\Re \z_n\to-\infty$, then $h(\z_n)\to p_{k+1}$.

In order to prove that $p_k\neq p_{k+1}$ and that $p_k$ is monotone, we notice that the length of the curve $h(z(\z))$ is infinite where $\z=u+\i v_1$ with $v_1$ to be a constant and $-\infty<u<\infty$. Moreover, by \cite{Wf, M} or \cite{HTTW, p.109}, the geodesic curvature of this curve is bounded by $\epsilon$ provided $v_1$ is large. From this, it is easy to see that $p_k\neq p_{k+1}$. Since $h$ is an orientation preserving diffeomorphism, we conclude that $p_k\neq p_j$ if $k\neq j$, and $p_k$ is monotone on $\Bbb S^1$. In particular, $p_+=\lim_{k\to\infty}p_k$ and $p_-=\lim_{k\to-\infty}p_k$
exist.

To prove (ii), we observe that for any $C>0$ there is $v_0>0$ independent of $k$ such that the $\Phi$-radius of $\z\in \Cal R_1$ is larger than $C$ for all $\z$ with $\Im \z>v_0$. By the argument in \cite{HTTW, p.102}, we conclude that for any  $\epsilon>0$, there is $v_0>0$ independent of $k$ such that if $\Im \z>v_0$, then $d_{\H^2}(h(\z)),\gamma_k)\le \epsilon$, where $\gamma_k$ is the geodesic joining $p_k$ and $p_{k+1}$. From the proof of Lemma 1.2, we see that given $v_0$, there exists $a>0$ independent of $k$ such that if $z\in \Cal S_2$ and $\Re z>a$, then $\Im \z(z)>v_0$.  From this we can conclude that (ii) is true.

In order to prove (iii), let $\delta>0$ as above but small and let $b=a_2$ which is in the definition of $\Cal S_2$. Suppose   $z_n\in\C_\alpha$ with $\Re z_n>b$. Let $k_n$ be such that $2k_n\pi\le \Im z_n<2(k_n+1)\pi$. Then $\lim_{n\to\infty}k_n=\infty$. For each $n$, let $\z=\z^{(k_n)}$ as above then $\z_n=\z(z_n)$ can be defined and $\z_n\in \Cal R_2$. By Lemma 1.1, for all $\z\in \Cal R_2$ with $\Im \z>0$ and $\Re\z=\Re\z_n$, $d_{\H^2}(h(\z_n),h(\z))\le C_2$ for some constant $C_2$ independent of $n$. From (ii), we conclude that $d_{\H^2}(h(\z_n),\gamma_{k_n})\le C_3$ for some constant $C_3$ independent of $n$. From this, the result follows.
\enddemo

\proclaim{Theorem 3.1} Let $h$ be an orientation preserving harmonic diffeomorphism from $\C$ into $\H^2$ with Hopf differential $\Phi=\phi dz^2=P\e(Q)dz^2$ such that
$|\partial h|^2|dz|^2$ is complete on $\C$ and such that 
\roster
\item"{(i)}" $Q(z)=z^n+\sum_{j=1}^{n-1}a_jz^{n-j}$ is a polynomial of degree $n\ge 1$;
\item"{(ii)}" $P\not\equiv0$ is an  entire function with order $\rho<n$; and
\item"{(iii)}"  there exists $\frac{\pi}{2n}>\delta>0$ and $R_0>0$ such that 
$$
\Sigma\cap 
\left\{z|\ |z|>R_0 \ \text{and} \ |\arg z-\frac{2k\pi}{n}| <\frac{\pi}{2n}+\delta\right\}
=\emptyset 
$$
for all $0\le k\le n-1$, where  $\Sigma$ is the set of all zeros of $P$.
\endroster 
Then the closure of the image of $h$ is the convex hull of a countable set $A$ of $\partial \H^2$ with exactly $n$ accumulation points.
\endproclaim
\demo{Proof} We claim that for any $\epsilon>0$ with $n\epsilon<\frac\pi2$, there exists a constant $C_1>0$ such that the maximal $\Phi$-radius $R_z$ of $z$ satisfies
$$
R_z\le C_1\tag3.1
$$ 
for all $z\in W_k$, $0\le k\le n-1$, where $W_k$ is the wedge
$$
W_k=\lf\{z\big|\ \lf|\arg z-\frac{(2k+1)\pi}{n}\ri|<\frac\pi{2n}-\epsilon\ri\}.
$$
for $0\le k\le n-1$. To prove the claim, note that there exists $\tau>0$ such that for   $z\in W_k$,
$\Re (z^n)\le -\tau |z|^n$. By the assumptions (i) and (ii), for any $z\in W_k$, let $\gamma $ be the half ray $\gamma(t)=t\e(\i \arg z)$ for $t\ge |z|$, then 
$$
\split
\int_{|z|}^\infty|\phi|^\frac12(\gamma(t))dt&\le \int_{|z|}^\infty \e\lf(-\frac12\tau t^n+C_2(1+t^{n-1}+t^{\tilde{\rho}})\ri)dt, \text{ for any } \rho<\tilde{\rho}<n,\\
&\le C_3
\endsplit
$$
where $C_2$ and $C_3$ are constants independent of $z$. Hence the maximal $\Phi$-radius of $z\in W_k$ is uniformly bounded. This proves the claim.

Next, for  each $0\le k\le n-1$, and for $\delta>4\epsilon>0$, let 
$$
V_{k,\epsilon}=\lf\{z\big|\ \lf|\arg z-\frac{2k\pi}{n}\ri|<\frac\pi{2n}+\epsilon\ri\}.
$$
 Define $V_{k,4\epsilon}$ similarly. By assumption (iii),  we can take a branch of $\log P$ in $\{z\in V_{k,4\epsilon}|\ |z|>R_0\}$. As in the proof of Theorem 1.1, there exist positive constants $R_2> R_1>R_0 $, $T_2>T_1$, $\epsilon_1>0$ and a conformal map $\z_k(z)$ which is of polynomial growth as a function of $z$ and   which will map $\Cal S^{(k)}_1=\{z\in V_{k,4\epsilon}|\ |z|>R_1\}$ injectively onto its image. For simplicity, we write $\z=\z_k$.  Moreover, if
$$\Cal S^{(k)}_2=\{z\in V_{k,\epsilon}|\ |z|>R_2\}
$$
$$\Cal R_1=\{\z|\ |\arg\z|<\frac\pi 2+2\epsilon_1\ \text{and}\ |\z|>T_1\}
$$
and
$$\Cal R_2=\{\z|\ |\arg\z|<\frac\pi 2+ \epsilon_1\ \text{and}\ |\z|>T_2\}$$
then $\z(\Cal S^{(k)}_1)\supset \Cal R_1 \supset \z(\Cal S^{(k)}_2)\supset \Cal R_2$.  Moreover, in $\Cal R_1$ the Hopf differential of $h$ is of the form $\Phi=\e(\z+Q_1(\z))d\z^2$ where
$Q_1(\z)\to0$ as $\z\to\infty$. Choose $a>b>T_2$.  As in the proof of (3.1), we have 
$$
R_z\le C_2\tag3.2
$$
for some constant $C_2$ for all $z\in \Cal S^{(k)}_2\cap \z^{-1}(\{\Re \z\le a\})$. Moreover, on $\Re \z=b$, $|\e(\z+Q_1(\z))|\ge C_3$ for some positive constant $C_3$. Hence if $\tilde w=\log|\partial_{\z}h|$ and if $e^{2\tilde v}|d\z|^2$ is the hyperbolic metric on $\z(\Cal S^{(k)})\cap \{\Re\z>a\}$, then $\tilde w\ge \tilde v-C_4$ for some constant $C_4$ because $e^{-2\tilde w(\z)}|\e(\z+Q_1(\z))|<1$ and $\tilde v\le C$ on $\Re z=b$ for some positive constant $C$. Let $\Gamma_k=\z^{-1}(\{\Re z=a\})$ and $\gamma_k=\z^{-1}(\{\Re \z=b\})$. Note that for fixed $c>T_2$, $\arg(\z^{-1}(c+\i t))\to (2k\pm\frac12) \pi/n$ as $t\to\pm\infty$. Let $\Omega$ be the component containing the origin of $\C\setminus\cup_{k=0}^{n-1}\Gamma_k$,  and let $\Omega_1$ be the component containing the origin of $\C\setminus\cup_{k=0}^{n-1}\gamma_k$. By (3.1) and (3.2) if we choose $\epsilon>0$ in (3.2) and then choose $\epsilon>0$ in (3.1) small enough then we have $R_z\le C_1+C_2$ for all $z\in \Omega$, and if $e^{2v}|dz|^2$ is the hyperbolic metric on $\Omega$ then $w=\log|\partial_z h|\ge v-C$ for some constant $C$ for all $z\in \partial \Omega_1$. Here we have used the fact that the hyperbolic metric on $\Omega$ is dominated by the hyperbolic metric on its subdomain. 

Next we want to show that $\inf_{z\in\C\setminus \Omega_1}R_z>0$. In fact, if $z\in \C\setminus\Omega_1$, then there is $k$ such that $\Re\z_k(z)\ge b$. Apply Corollary 2.1 on the disk with center $\z_k$ and radius 1, we can conclude that   on $\Re \z_k\ge b$ the maximal $\Phi$-radius is bounded below by a positive constant   independent of $\z_k$, because $\Phi=\e(\z_k+o(1))d\z_k^2$. 

 Since  $\Omega$ is strongly hyperbolic and $|\partial h|^2|dz|^2$ is complete in $\C$,  $h$ is quasi-conformal on $\Omega_1$ by Theorem 2.1 and Remark 2.1.

On the other hand, by Lemma 3.1, if we choose $a$ and $b$ large enough, then for each $k$,  there exist $p_j^{(k)}\in \partial \H^2$, $j\in \Bbb Z$, which are monotone in $\Bbb S^1$  such that the intersection of the closure of the image under $h$ of the set $\{\z\in\Cal R_2,\, \Re\z\ge b\}$  with $\partial \H^2$ is equal to
$$
\Cal A=\bigcup_{k=0}^{n-1}\{p_j^{(k)}|\ j\in\Bbb Z\}\cup\bigcup_{k=0}^{n-1}\{p_+^{(k)}, p_-^{(k)}\},
$$
 where $p_\pm ^{(k)}=\lim_{j\to\pm\infty}p_j^{(k)}$. Moreover, if $\z_n\in \Cal R_2$ with $\Re \z_n\ge b$ and $\Im z_n\to+\infty$ (respectively $\Im z_n\to-\infty$) then $h(\z_n)\to p_+^{(k)}$ (respectively $h(\z_n)\to p_-^{(k)}$).

Since $h$ is at most linear growth in $\Cal R_2$ with respect to $\z$, $h$ is of polynomial growth on $V_{k,\epsilon}$, provided $\epsilon>0$ is small enough. It is easy to see that $h$ is at most of linear growth on $W_k$. By the definition of $\Omega_1$, we see that $h$ is of polynomial growth on $\Omega_1$. Namely, there exist positive constants  $\ell$ and $C$ such that 
$$
d_{\H^2}\lf(h(z),o\ri)\le C\lf(d_{\C}(z,0)+1\ri)^\ell\tag3.3
$$
for all $z\in \Omega_1$, where $o$ is a fixed point in $\H^2$ and $0$ is the origin of $\C$. We claim that the image of $h$ is the convex hull of 
$\Cal A$ together with at most finitely many points $q_j\in \partial\H^2$. By theorem 4.8 in  \cite{C-T} and theorem 5 in \cite{Wn}, it is sufficient to show that $\overline{h(\C)}\cap \partial\H^2$ is $\{p_j^{(k)}|\ j\in\Bbb Z\}\cup\{p_+^{(k)}, p_-^{(k)}\}$ together with at most finitely many points $q_j$. Suppose $q_1,\dots,q_m$ are distinct points in 
$$
\lf(\overline{h(\C)}\cap \partial\H^2\ri)\setminus \Cal A.
$$
 There exist disjoint neighborhoods $U_1,\dots,U_m$ of $q_1,\dots,q_m$ respectively in $\overline\H^2$. We may choose $U_j$, $1\le j\le m$ small enough so that $h^{-1}(U_j)\subset \Omega_1$. For if this is not true, then there exists   $q_j$ and a sequence of neighborhoods $U_{j,n}$ such that $\bigcap_{n=1}^\infty U_{j,n}=\{q_j\}$ and such that $h^{-1}(U_{j,n})$ is not contained in $\Omega_1$ for each $n$. By choosing a subsequence, we may assume that there is $z_n\in U_{j,n}$ such that $\Re\zeta_k(z_n)\ge b$ under the map $\z_k$ described above. Since $h(z_n)\to q_j$ by construction, we conclude that $q_j$ must be $p_l^{(k)}$ or $p^{k}_\pm$ for some $k$ and $l$. This is a contradiction. Hence we may choose $U_j$ such that $h^{-1}(U_j)$ is contained in $\Omega_1$. Moreover, we may assume that $U_j$ is bounded by a geodesic line in $\H^2$. Let  $f_j(z)=d_{\H^2}(u(z), \H^2\setminus U_j)$, then $f_j$ harmonic because $d_{\H^2}(\cdot, \H^2\setminus U_j)$  is convex by \cite{B-O}. Note that $f_j$ is smooth in $h^{-1}(U_j)$,    $f_j(z)=0$ for $z\in \C\setminus h^{-1}(U_j)$ and there exists a constant $C_3$
$$
f_j(z)\le C_3\lf(d_\C(z)+1\ri)^\ell
$$
  for all $z$ and for all $1\le j\le m$ by (3.3). Since $h^{-1}(U_j)$, $1\le j\le m$, are disjoint and    nonempty, $m$ is bounded from above by a constant depending only on $\ell$ by Theorem 3.4 in \cite{L-W 1}. This proves the claim. 

Observe that  each $q_j$ must lie between $p_+^{(k)}$ and $p_-^{(k+1)}$ for some $k$.  Here we use the convention that $p_+^{(n)}=p_-^{(0)}$.  
Since if $ \gamma_k(t) =\z^{-1}(b+\i t)$, then $\lim_{t\to\pm\infty}h(\gamma_k(t))= p_{k,\pm}$, we conclude that $h(\Omega_1)$ is bounded by $h(\gamma_k)$ and the geodesics  joining consecutive points of $p_+^{k}$,  $q_j$ and $p_-^{(k+1)}$, with $q_j$ between $p_+^{k}$ and $p_-^{(k+1)}$, and they are oriented positively. 
 Since $h$ is quasi-conformal on $\Omega_1$, 
for each $k$ if $\z_n\in \Cal R_1$ with $\Re\z_n\le b$ and $\Im z_n\to+\infty$ (respectively $\Im z_n\to-\infty$) then $h(\z_n)\to p_+^{(k)}$ (respectively $h(\z_n)\to p_-^{(k)}$). Again, using the fact that $h$ is quasi-conformal on $W_k$, we conclude that for $z\in W_k$, and if $z\to\infty$ then $h(z)$ will converge to a point $q_k$ in $\overline{\H^2}$. But $q_k$ must be equal to $p_+^{(k)}$ and $p_-^{(k+1)}$ at the same time. Hence the closure of the image of $h$ is  the convex hull of the set $A$ consisting of $p_j^{(k)}$, $q_k $ which is countable and has exactly $n$ accumulation points. 

\enddemo
It is clear that the theorem is true for any polynomial $Q$ without requiring the leading coefficient to be $1$ as long as the zeros of the entire function $P$ are distributed in the corresponding sections. For instance, we conclude immediately from the theorem the following.
 
\proclaim{Corollary 3.1} Let $h$ be an orientation preserving harmonic diffeomorphism from $\C$ into $\H^2$ with Hopf differential $\Phi=P\e(Q)dz^2$ where $P$ and $Q$ are polynomials with $\deg Q=n$. Suppose $|\partial h|^2|dz|^2$ is complete in $\C$. Then the image of $h$ is the convex hull of a countable set $A$ of $\partial \H^2$ with exactly $n$ accumulation points.
\endproclaim
 
Figures 1, 5, and 6 show the horizontal trajectories structures of holomorphic quadratic differentials which are included in the corollary 3.1. Figure 1 also shows the image of the  harmonic map corresponding to $e^zdz^2$ which is the basis of all the discussion in this paper.

\subheading{\S4 Images of harmonic diffeomorphisms with Hopf differential $f(e^z)dz^2$}

In this section, we will study the images of certain harmonic diffeomorphisms from $\C$ into $\H^2$ with Hopf differentials of the form $f(e^z)dz^2$. As before for $\alpha\ge 0$, let $\C_\alpha=\{ z=x+\i y\, |\, x>\alpha\}$.
\proclaim{Lemma 4.1} Let    $h: \C_0\to \H^2$ be an orientation preserving harmonic diffeomorphic injection with Hopf differential $\Phi=dz^2$. Suppose that for some $x_0>0$,  $\lim_{y\to+\infty}h(x_0+\i y)=p_1$ and $\lim_{y\to-\infty}h(x_0+\i y)=p_2$ for some $p_1$, $p_2$ in $\partial \H^2$. Then $p_1=p_2=p$, and for all $x_0>0$
 $$\lim_{{|z|\to\infty}\atop{\Re z\ge x_0}}h(z)=p.$$
\endproclaim
\demo{Proof} By proposition 1.5 in \cite{ T-W 1}, see also \cite{Wn}, the energy density of $h$ in the half-plane  $0<\Re z<\infty$ is bounded. Since for some $x_0>0$,  $\lim_{y\to+\infty}h(x_0+\i y)=p_1$ and $\lim_{y\to-\infty}h(x_0+\i y)=p_2$, where $p_1,\ p_2\in\partial \H^2$, we conclude that for all $x_1>x_0>0$, 
$$\lim_{{\Im z\to+\infty}\atop{x_0<\Re z<x_1}}h(z)=p_1$$
 and 
$$\lim_{{\Im z\to-\infty}\atop{x_0<\Re z<x_1}}h(z)=p_2.$$
Identify $\overline{\H^2}=\H^2\cup\partial\H^2$ with the unit disk.  We claim that for  any  $x_0>0$, the closure of $h(\C_{x_0})$ in $\overline{\H^2}$ is $\overline\Omega$ where $\Omega$ is  the domain bounded by the curve $h(x_0+\i y)$, $-\infty<y<\infty$ and one of the arc on $\Bbb S^1$ with end points $p_1$ and $p_2$. Obviously, $h(\C_{x_0})$ is contained in such an $\Omega$ because $h$ is injective. Suppose that the claim is not true, then there is $q$ on the boundary of $h(\C_{x_0})$ such that $q\in\H^2$ and there is a geodesic arc $\gamma$ in $h(\C_{x_0})\subset\H^2$ from a point $q_1$ in $h(\C_{x_0})$ to $q$ with $\gamma(\ell)=q$, where $\ell$ is the length of $\gamma$ and is finite. Without loss of generality, we may asumme that $\gamma([0,\ell))\subset h(\C_{x_0})$. Let $\beta=h^{-1}(\gamma)$. Then $\beta$ is a path in $\C_{x_0}$ such that $\beta(t)\to\infty$ as $t\to\ell$ because $q$ is in the boundary of $h(\C_{x_0})$. Moreover, $\Re \beta(t)\to +\infty$. Otherwise, we would have $\beta(t)\to p_1$ or $p_2$. However, the pull-back metric under $h$ is given by
 $ (e+2)dx^2+(e-2)dy^2$, where $e$ is the energy density of $h$, and $e>2$. We then have
$$\split
\ell&=\int_{0}^\ell\lf[ (e+2)\lf( \frac{dx}{dt}\ri)^2+(e-2)\lf( \frac{dy}{dt}\ri)^2\ri]^\frac12 dt\\
&\ge \sqrt 2\lf(\lim_{t\to\ell}x(\ell)-x(0)\ri)\\
&=\infty
\endsplit
$$
which is a contradiction. Hence $h(\C_{x_0})=\Omega$. Suppose $p_1\neq p_2$, then $\overline{h(\C_{x_0})}$ contains a nontrivial arc on $\Bbb S^1$. However,   for $x_0>0$, $h$ is of at most linear growth. By theorem 3.4 in \cite{L-W 1}, we conclude that 
$\overline{h(\C_{x_0})}\cap \partial\H^2$ consists of only finitely many points. Hence we must have $p_1=p_2=p$.  Since $h$ is a diffeomorphism, we must have $$\lim_{{|z|\to\infty}\atop{\Re z\ge x_0}}h(z)=p.$$
\enddemo

\proclaim{Lemma 4.2} Let $0<\beta\le\pi$ and let  $h:e^{\i\beta }\C_0 \to \H^2$ be an orientation preserving harmonic diffeomorphic injection with Hopf differential $\Phi=  dz^2$. Suppose that for some $x_0>0$, 
$$
\lim_{y\to\infty}h\lf(e^{\i\beta }(x_0+\i y)\ri)=p_1
\quad\text{and}\quad
\lim_{y\to-\infty}h\lf(e^{\i\beta}(x_0+\i y)\ri)=p_2
$$
for some $p_1$, $p_2$ in $\partial \H^2$. Then $p_1\neq p_2$ and  for all $x_0>0$,  $\overline{h(e^{\i\beta}\C_{x_0})}\cap\partial\H^2=\{p_1,p_2\}$.
\endproclaim
\demo{Proof} As in the proof of Lemma 4.1, we conclude that for any $x_0>0$, 
$$\lim_{y\to\infty}h\lf(e^{\i\beta}(x_0+\i y)\ri)=p_1,$$
 and 
$$\lim_{y\to-\infty}h\lf(e^{\i\beta}(x_0+\i y)\ri)=p_2.$$ Let $x_0>0$. Since $0<\beta\le \pi$, by Lemma 1.1, suppose $z_n\in e^{\i\beta}\C_{x_0}$, if $\Re z_n\to-\infty$,  then $\lim_{n\to\infty}h(z_n)=p_1$; and if  $\Re z_n\to\infty$,  then $\lim_{n\to\infty}h(z_n)=p_2$. If $z_n\to\infty$ and for all $n$, $x_0\le \Re z_n<x_1$   for some $x_1$, then $h(z_n)$ are uniformly bounded. Hence $\overline{h(e^{\i\beta}\C_{x_0})}\cap\partial_{\infty}\H^2=\{p_1,p_2\}$. 

To prove that $p_1\neq p_2$. Note that $\i n\in e^{\i\beta}$ for any positive integer $n$.    Moreover, it is easy to see that $z_n=e^{\i\frac \beta3}\i n=-n\sin\frac\beta3+\i n\cos\frac\beta3$ and $\tilde z_n=n+\i n\cos\frac\beta3$ are  in $e^{\i\beta}C_{x_0}$ if $n$ is large. Let $L_n$ be the horizontal line joining $z_n$ and $\tilde z_n$. By the arguments in section 3 of \cite{HTTW}, we conclude that $h(L_n)$ is of uniformly bounded distance from the geodesic passing through $h(z_n)$ and $h(\tilde z_n)$. Since $h(z_n)\to p_1$ and $h(\tilde z_n)\to p_2$, if $p_1=p_2=p$ then $h(L_n)\to p$ as $n\to\infty$.  
On the other hand,    $h(\i n\cos\frac\beta3)$ are uniformly bounded. This is a contradiction. Therefore, $p_1\neq p_2$.\enddemo

\proclaim{Theorem 4.1}
Let $m$, $n$ be nonnegative intergers and let $P(t)$ be a nonconstant rational function of the form 
$$
P(t)=\sum_{k=-m}^{n} a_kt^k,
$$
with $a_{-m}\neq0\neq a_n$. Suppose   $h$ is an orientation preserving harmonic diffeomorphism from $\C$ into $\H^2$ with Hopf differential given by
$$
\Phi=P(e^z)dz^2 
$$
such that $|\partial h|^2|dz|^2$ is a complete metric. Then   $\Cal A=\overline{h(\C)}\cap\partial\H^2$ is countable which has exactly one accumulation point if $m$ or $n=0$, and  $a_0\ge0$; and has two accumulation points otherwise.  Moreover $\overline{h(\C)}$  is the convex hull of $\Cal A$.
\endproclaim
\proclaim{Remark}
Figures 1 and 2 in the appendix show horizontal trajectories structures for the case that $m$ or $n=0$, and  $a_0\ge0$. In fact, in both figures, $m=0$, and $a_0=0$ and $1$ respectively. The other case are showed by the figures 3 and 4. In figure 3, $m=0$ but $a_0=-1$. The image of the corresponding harmonic map has 2 accumulations both are limits from 1 side. In figure 4, both $m$ and $n$ are not zero and the image of the corresponding harmonic map has 2 accumulations both are limits from 2 sides.
\endproclaim
\demo{Proof} Suppose that $m>0$ and $n>0$. By the proof of Lemma 3.1, we can conclude that
there exist  $p_k$, $k\in \Bbb Z$ such that
$\overline{h(\{z|\   \Re z\ge0\}}\cap \partial \H^2$ is equal to $\overline{\{p_k\}_{k\in \Bbb Z}}$ and the $p_k$ are monotone on $\Bbb S^1$. Moreover, if $p_k\to p_\pm$ as $k\to\pm\infty$, then $\lim_{{\Im z\to +\infty,\ \Re z\ge0}}h(z)=p_+$ and $\lim_{{\Im z\to -\infty,\ \Re z\ge0}}h(z)=p_-$. 

Similarly, there exist  $q_k$, $k\in \Bbb Z$ such that
$\overline{h(\{z|\   \Re z\le0\}}\cap \partial \H^2$ is equal to $\overline{\{q_k\}_{k\in \Bbb Z}}$ and the $q_k$ are monotone on $\Bbb S^1$, and if $q_k\to q_\pm$ as $k\to\pm\infty$, then $\lim_{{\Im z\to +\infty,\ \Re z\le0}}h(z)=q_+$ and $\lim_{{\Im z\to -\infty,\ \Re z\le0}}h(z)=q_-$. Hence $q_+=p_+$ and $q_-=p_-$. Since $h$ is a diffeomorphism, $p_+\neq p_-$ and $\Cal A=\{p_k,\ q_k\}_{k\in\Bbb Z}\cup\{p_+,p_-\}$, which has two accumulation points.

Next, let us consider the case that $m$ or $n=0$. Without loss of generality, we may assume that $m=0$.   As before, there exist  $p_k$, $k\in \Bbb Z$ such that
$\overline{h(\{z|\   \Re z\ge0\}}\cap \partial \H^2$ is equal to $\overline{\{p_k\}_{k\in \Bbb Z}}$ and the $p_k$ are monotone on $\Bbb S^1$. Let $p_+$ and $p_-$ defined as above.

 Suppose $a_0=0$. Then  we can conclude as in the proof of Theorem 3.1 that $p_+=p_-=p$ and $\Cal A=\{p_k\}_{k\in\Bbb Z}\cup\{p\}$ which has only one accumulation points. 

Suppose $m=0$ and $a_0\neq 0$, let $a_0=\rho^2 e^{2\i \beta}$ with $0\le \beta<\pi$, $\rho>0$. There exists $\delta>0$ such that if $|t|<\delta$, we can take a branch of the square root of $P(t)$ and 
$$
\sqrt{P(t)}= \rho e^{\i\beta}+tg(t),
$$
where $g(t)$ is analytic and 
$$|g(t)|\le C_1\tag4.1
$$ for some constant $C_1$ for $|t|<\delta$. Let $\tilde g(t)$ be such that $\tilde g'=g$ on $|t|<\delta$ and $\tilde g(0)=0$. Let $x_0<0$ be small enough so that $|e^z|<\delta$ on $\Re z\le x_0$. Define
$$
\split
\z(z)&=\int_{x_0}^z\sqrt{P(e^\xi)}d\xi\\
&=\rho e^{\i\beta}(z-x_0)+\int_{x_0}^ze^\xi g(e^\xi)d\xi\\
&=\rho e^{\i\beta} z+\tilde g(e^z)+\z_0
\endsplit
$$
for all $z$ with $\Re z\le x_0$, where $\z_0$ is a constant. Here the integration is along the straight line from $x_0$ to $z$.  Then $\z$ is analytic. By (4.1), if we choose $x_0$ small enough, then $\z$ is injective. Since $|\tilde g(e^z)|\le C_2|e^z|$ for some constant $C_2$, if we choose $x_0$ small enough, then the analytic map $z\to \z_1=-(\z-\z_0)$ will map $\{\Re z\le x_0\}$ injectively onto its image $\Cal R$. Moreover
$$
e^{\i\beta}\{\z_1|\ \Re z_1\ge \rho x_0+1\}\subset \Cal R\subset
  e^{\i\beta}\{\z_1|\ \Re z_1\ge \rho x_0-1\}.
$$
The Hopf differential of $h$ in the $\z_1$ plane is given by $d\z_1^2$. As before, we have $$\lim_{y\to\pm\infty}h(x_0+\i y)=p_\pm.$$ Hence if $\beta=0$, we have $p_+=p_-=p$ by Lemma 4.1 and $\Cal A$ is countable with only one accumulation point. If $\beta>0$, we have $p_+\neq p_-$ by Lemma 4.2 and  $\Cal A$ is countable with exactly two accumulation points. The last statement of the theorem follows from   theorem 4.8 in \cite{C-T} and theorem 5 in \cite{Wn}. 
\enddemo
As an application, we use Theorem 4.1 to study harmonic diffeomorphic injection from a flat cylinder to  a hyperbolic cylinder. Let $N$ be a hyperbolic cylinder and  Let $\C^*=\C\setminus\{0\}$.  Let $\Phi(\C^*,N)$ be the set of all Hopf differentials of orientation preserving harmonic diffeomorphic injections $h$ from $\C^*$ to $N$ such that $|\partial h|^2|dz|^2$ is complete on $\C^*$. Let $\Cal P(N)$ be the set of holomorphic quadratic differentials on $\C^*$ defined by 
$$
{\Cal P}(N)=\left\{\frac{P(z)}{z^2}dz^2\big|\ P(z)=\sum_{k=-m}^{n} a_kz^k\ \text{for some $0\le m,n\in\Z$, and $ P\neq a_0$}  \right\}.
$$

 $$\Cal P_1(N)=\left\{\frac{P(z)}{z^2}dz^2\in\Cal P(N)\big|\ P(z)=\sum_{k=-m}^{n} a_kz^k\ \text{with $m$ or $n=0$, and $a_0\ge0$}  \right\}
$$
and
 $\Cal P_2(N)=\Cal P(N)\setminus \Cal P_1(N)$.

\proclaim{Corollary 4.1} With the above notations we have $\Phi(\C^*,N)\cap \Cal P(N)$ is either a subset of $\Cal P_1(N)$ or a subset of $\Cal P_2(N)$. Moreover,  if $\Phi(\C^*,N)\cap \Cal P(N)\neq\emptyset$, then it is   a subset of $\Cal P_1(N)$ if and only if $N$ has a cusp.
\endproclaim
\demo{Proof} Let $z^{-2}P(z)dz^2\in \Phi(C^*,N)$ be the Hopf differential of an orientation preserving harmonic diffeomorphic injections $h$ from $\C^*$ into $N$.  Lifting $h$ to the universal coverings, we have an orientation preserving harmonic diffeomorphic injection, denoted by $h$ again, from $\C$ into $\H^2$, with Hopf differential given by
$$
P(e^z)dz^2
$$
and an element $\rho$ of the M\"obius group which generates $\pi_1(N)$ such that
$$
h(z+2\pi i)=\rho(h(z)).
$$
Note that $|\partial h|^2|dz|^2$ is complete on $\C$. Let  $\Cal A=\overline{h(\C)}\cap\partial_{\infty}\H^2$.   Since $h$ is equivariant, $\Cal A$ is invariant under $\rho$. This implies that the set of  fixed points of  $\rho$ is exactly the set of accumulation points of $A$. The corollary then follows easily from Theorem 4.1.

\enddemo

\proclaim{Remark 4.1}
It was proved in \cite{Wn,W-A,T-W 1} that given a holomorphic quadratic differential $\Phi$ on $\C$ or on $\H^2$ there exists an orientation preserving harmonic diffeomorphic injection from $\C$ or $\H^2$ to $\H^2$ whose Hopf differential is the given $\Phi$. Corollary 4.1   shows  that the prescribed Hopf differential problem is not alway solvable from $\C^*$ into $N$ where $N$ is a  hyperbolic cylinder.
\endproclaim

Our next result is to consider the image of a harmonic map with Hopf differential with infinite order.

\proclaim{Theorem 4.2}
Let $h$ be an orientation preserving harmonic diffeomorphic injection from  $\C$ into $\H^2$ such that $|\partial h|^2|d z|^2$ is complete. Suppose that the Hopf differential of $h$ is given by
$$
\Phi=\e^{(k)}(z)dz^2,
$$
for some positive integer $k$, where $\e^{(k)}(z)$ is defined inductively by   $\e^{(0)}(z)=1$ and $\e^{(j)}(z)=\e(\e^{(j-1)}(z))$.
Let $\Cal A=\overline{h(\C)}\cap\partial\H^2$. Then $\Cal A=\cup_{j=0}^k\Cal A_j$ such that
\roster
\item  $A_j$ is countable and discrete for each  $0\le j\le k-1$;
\item   $\Cal A_j$ consists of all isolated accumulation points of $\Cal A_{j-1}$ for $1\le j\le k$;
\item $\Cal A_k$ consists of only one point.
\endroster
\endproclaim

Please see figure 7 in the appendix for the horizontal trajectories structure of $e^{e^z}dz^2$ and the corresponding image of the harmonic map.

\noindent
\demo{Proof} We may assume that $k\ge 2$ because   $k=1$ is a special case of Theorem 3.1.
First of all, we want to find out the domains such that $\Phi$ can be written in the form of Lemma 3.1.
Given any $\alpha\in\R$, $1\le l\le k-1$ and $(n_1, n_2,\ldots,n_l)\in\Z^l$, we define the open subsets $\cS_{(n_1,\ldots,n_l)}$ inductively by
$$
\cS_{(n_1)} =\{ z\in\C\,|\, \Re z>\alpha,\, |\Im z-2n_1\pi|<\pi\}
$$
and
$$
\cS_{(n_1,\ldots,n_l)}=\{ z\in\cS_{(n_1,\ldots,n_{l-1})}\, |\, \Re\zeta_{l-1}>\e^{(l-1)}(\alpha),\, |\Im\zeta_{l-1}-2n_l\pi|<\pi\},
$$
where $\zeta_{l-1}=\e^{(l-1)}(z)$. Then
 $\zeta_l=e^{\zeta_{l-1}}=\e^{(l)}(z)$ maps $\cS_{(n_1,\ldots,n_l)}$ one-one onto the open set
$$
\Omega_l=\C\setminus \left(\{\z_l\in\R| \z_l\le 0\}\cup \{\z_l|\ |\z_l|\le \e^{(l)}(\alpha)\}\right),
$$
and in terms of $\zeta_l$
$$
\Phi=\frac{\e^{(k-l)}(\zeta_l)}{\prod_{j=0}^{l-1}(\log^{(j)}\zeta_l)^2}  d\zeta_l^2.\tag4.2
$$
In particular, for $l=k-1$,
$$
\split
\Phi&=\frac{\e(\zeta_{k-1})}{\prod_{j=0}^{k-2}(\log^{(j)}\zeta_{k-1})^2}  d\zeta_{k-1}^2\\
 &=\e\lf(\zeta_{k-1}-2\sum_{j=1}^{k-1}\log^{(j)}\zeta_{k-1}\ri) d\zeta_{k-1}^2
\endsplit\tag4.3
$$
on
$$
\Omega_{k-1}=\C\setminus \left(\{\z_{k-1}\in\R| \z_{k-1}\le 0\}\cup \{\z_{k-1}|\ |\z_{k-1}|\le \e^{(k-1)}(\alpha)\}\right).
$$
A further transformation 
$$
\eta=\zeta_{k-1}-2\sum_{j=1}^{k-1}\log^{(j)}\zeta_{k-1}
$$
will put the Hopf differential into the form of Lemma 3.1 and we can conclude on the boundary behaviour of the harmonic map $h$. However, to ensure that there are no other ideal boundary point, we need to show that $h$ is quasiconformal in certain domain.

In order to do so, given any $\beta>>1$, we define $E_{\beta}=\{ \z_{k-1}\in\Omega_{k-1}\, |\, \Re\eta(\z_{k-1}) > \beta\}$ and claim that for any $\alpha\in\R$, there are simply-connected domains
$V_0\subset\tilde{V}_0\subset\{ z\in\C\, |\, \Re z>\alpha -1\}$ such that
$$
\C\setminus \tilde{V}_0=\{ \Re z<\alpha-1\}\cup \left(\bigcup_{(n_1\ldots,n_{k-1})} \tilde{T}_{(n_1\ldots,n_{k-1})}\right)
$$
and
$$
\C\setminus V_0=\{ \Re z<\alpha\}\cup \left(\bigcup_{(n_1\ldots,n_{k-1})} T_{(n_1\ldots,n_{k-1})}\right),
$$
where $\tilde{T}_{(n_1\ldots,n_{k-1})}$, respectively $T_{(n_1\ldots,n_{k-1})}$, is the component of the preimage of $E_{\beta+1}$, respectively $E_{\beta}$,  under the map $\e^{(k-1)}(z)$ corresponding to the branch of $\log$ given by $(n_1,\ldots,n_{k-1})$. Moreover, there are constants $C_0$, $M_0$, $\delta $ with $C_0$ and $\delta>0$ such that
$$
\sup_{\tilde{V}_{0}}R_z\le C_{0},\quad 
\inf_{\C\setminus V_{0}}R_z\ge \delta,\quad\text{and}\quad
\inf_{\partial V_{0}}(w-v)(z)\ge M_0, \tag 4.4
$$
where $R_z$ is the maximal $\Phi$-radius at $z$, $w=\log|\partial h|$ and $e^{2v}|dz|^2$ is the Poincar\'e metric on $\tilde{V}_0$.
If the claim is true, then the last inequality of (4.4) implies that $w\ge v-M_0$ for all $z\in V_0$, and hence, by Theorem 2.1, one can conclude that $h$ is quasiconformal on $V_0$.

To prove the claim, we note that for $z\in T_{(n_1,\dots,n_{k-1})}$ then the image of $z$ under $\exp^{(k-1)}$ is in  $E_\beta$. By Corollary 2.1 as in the proof of Theorem 3.1, we conclude that $R_z\ge \delta>0$ for some $\delta>0$ independent of $(n_1,\dots,n_{k-1})$. On the other hand, since $e^{(k)}(z)\to e^{(k-1)}(0)$ if $\Re z\to-\infty$, we also have $R_z\ge \delta$ by Corollary 2.1 if $\Re z\le \alpha$ by choosing a possible smaller $\delta$. The second inequality of (4.4) is proved.

Let $z\in \partial V_0$, then either $\Re z=\alpha$ or the image of $z$ under $\exp^{(k-1)}$ is on the  boundary of $E_\beta$. In the first case, $e^{2w(z)}\ge |e^{(k-1)}(z)|\ge C$ for some constant $C>0$ independent of $z$. Hence it is easy to see that $w(z)-v(z)\ge M_0$ for some constant $M_0$ because $\tilde V_0$ is strongly hyperbolic. In the second case, then we can proceed as in the proof of Theorem 3.1 and obtain the third  inequality in (4.4).

To prove the first inequality in (4.4), 
 we let
$$
V_{k-1}=\Omega_{k-1}\setminus E_{\beta}\subset \tilde{V}_{k-1}=\Omega_{k-1}\setminus E_{\beta+1}.
$$
Then it is easy to see that $V_{k-1}\subset \tilde{V}_{k-1}$ are simply-connected domains in $\Omega_{k-1}$ and there is $C_{k-1}$ such that
$$
\sup_{\tilde{V}_{k-1}}R_{z}\le C_{k-1}.\tag4.5
$$
In fact, for all $z\in\tilde{V}_{k-1}$, there is a divergent path $\gamma$ in $\tilde{V}_{k-1}$ such that
$$
L_{\Phi}(\gamma)< C_{k-1}.
$$

Now, for $l=k-2$, we consider subsets in $\Omega_{k-2}$ containing the preimage of $V_{k-1}$ and $\tilde{V}_{k-1}$ under the exponential map $\zeta_{k-1}=\exp(\zeta_{k-2})$. It is clear from the property of the exponential map that
$$
\split
{V}_{k-2}&=\e^{-1}({V}_{k-1})\cup\left[ \left(\ol{V_{k-2}^+}\cup \ol{V_{k-2}^-} \right)\cap\Omega_{k-2}\right]\\
\tilde{V}_{k-2}&=\e^{-1}(\tilde{V}_{k-1})\cup\left[ \left(\ol{V_{k-2}^+}\cup \ol{V_{k-2}^-} \right)\cap\Omega_{k-2}\right]
\endsplit
$$
are simply-connected domains in $\Omega_{k-2}$ such that
$$
V_{k-2}^{\pm}\subset V_{k-2}\subset \tilde{V}_{k-2},
$$
where, for $l=1,\ldots,k-1$,
$$
\split
V_{l-1}^+ &=\{ z\in\cS_{(n_1,\ldots,n_{l-1})}\, |\, \Re\zeta_{l-1}<\e^{(l-1)}(\alpha),\, \Im\zeta_{l-1}>0\} \\
V_{l-1}^- &=\{ z\in\cS_{(n_1,\ldots,n_{l-1})}\, |\, \Re\zeta_{l-1}<\e^{(l-1)}(\alpha),\, \Im\zeta_{l-1}<0\}.
\endsplit
$$
We note that, for $l=1,\ldots,k-1$,
$$
\cS_{(n_1,\ldots,n_{l-1})}=V_{l-1}^+\cup V_{l-1}^- \cup \left[ \cup_{n_l\in\Z}\left( \overline{\cS_{(n_1,\ldots,n_{l-1},n_l)}}\cap\cS_{(n_1,\ldots,n_{l-1})} \right) \right].
$$
We want to show that there exists $C'_{k-2}>0$ such that for all $z\in\tilde{V}_{k-2}$, there is a divergent path $\gamma$ in $\tilde{V}_{k-2}$ with $L_{\Phi}(\gamma)<C'_{k-2}$. This will
immediately implies that
$$
\sup_{\tilde{V}_{k-2}}R_z\le C_{k-2}
$$
for some $C_{k-2}>0$.
To prove this, we note that for all $r_0>\e^{(k-2)}(\alpha)$,
$$
\split
\int_{r=r_0,\, 0<\theta<\pi}|\Phi||d\zeta_{k-2}|&\le 
C\int_{r=r_0,\, 0<\theta<\pi} \frac{rd\theta}{|\zeta_{k-2}||\log\zeta_{k-2}|\cdots|\log^{(k-2)}\zeta_{k-2}|}\\
&\le \frac{C}{\log r_0 \cdots\log^{(k-2)}r_0}\\
&\to 0 \hbox{ as } r_0\to+\infty .
\endsplit
$$
Then for all point $\zeta_{k-2}\in V_{k-2}^{\pm}$, it can be connected to a point on the vertical line $\{ \Re\zeta_{k-2}=\e^{(k-2)}(\alpha)\}$ by a circular arc with uniformly bounded $\Phi$-length. By using $\zeta_{k-1}=\e({\zeta_{k-2}})$ to map a point on $\partial\Omega_{k-1}$, we can find a divergent path in $\tilde{V}_{k-1}$ with $\Phi$-length bounded by $C_{k-1}$. Lifting this path to $\tilde{V}_{k-2}$ and together with the circular arc, we find a path starting from any point in $V_{k-2}^{\pm}$ a divergent path in $\tilde{V}_{k-2}$. The same is obviously true for other $z\in \tilde{V}_{k-2}$ since they   belong to $\e^{-1}(\tilde{V}_{k-1})$. This proves our assertion that
$$
\sup_{\tilde{V}_{k-2}}R_z\le C_{k-2}
$$
for some $C_{k-2}>0$.

Continue in this way,  for all $l\le k-1$  we can define $V_l$, $\tilde V_l$  and $V_l^{\pm}$, such that
$$
V_0=\e^{-1}(V_1)\quad\text{and } \tilde{V}_0=\e^{-1}(\tilde{V}_1)\cup \{z\,|\, \alpha-1<\Re z<\alpha\}.
$$
Moreover, we can prove inductively that 
$$
\sup_{\tilde{V}_{l}}R_z\le C_{l} 
$$
for some constant $C_l$. Finally, it is easy to see that the $R_z\le C$ for some contant $C $ if $\alpha-1<\Re z<\alpha$ because such a point can be joined by a line with bounded $|\phi|$-length to a point in $\e^{-1}(\tilde V_1)$. This completes the proof of the first inequality in (4.4).

Now we can study the structure of the boundary points of $h(\C)$. Firstly, for any $(n_1,\dots,n_{k-1})\in\Z^{k-1}$, using  Lemma 3.1, we can argue as before to conclude that there exists 
monotone sequence $p_{(n_1,\dots,n_{k-1});j_0}$, $j_0\in \Z$ and $p_{(n_1,\dots,n_{k-1});+}$, $p_{(n_1,\dots,n_{k-1});-}$ in $\partial \H^2$ such that, for $\beta$ sufficiently large,
$$
\split
&\ol{h(\cS_{(n_1,\dots,n_{k-1})}\cap\{z|\ \Re\eta(z)\ge \beta\})}\cap\partial\H^2\\
&\qquad=\{p_{(n_1,\dots,n_{k-1});j_0}\}_{j_0\in\Z}\cup\{p_{(n_1,\dots,n_{k-1});+},p_{(n_1,\dots,n_{k-1});-}\}
\endsplit\tag4.6
$$
$$
\lim_{j_0\to\pm\infty}p_{(n_1,\dots,n_{k-1});j_0}=p_{(n_1,\dots,n_{k-1});\pm}.
$$
Moreover, if $z_n\in\cS_{(n_1,\dots,n_{k-1})}\cap\{z|\ \Re\eta(z)\ge \beta\}$ and $z_n\to\infty$ then
$$
h(z_n)\to
\cases
p_{(n_1,\dots,n_{k-1});+},& \text{if  }\Im\eta(z_n)\to\infty\\
p_{(n_1,\dots,n_{k-1});-},& \text{if  }\Im\eta(z_n)\to-\infty
\endcases
$$ 
In fact,  we conclude by (4.3) that the energy density of $h$ is bounded on the set $z\in \cS_{(n_1,\dots,n_{k-1})}$ such that $a<\Re\eta(z)<b$ and $|\Im\eta(z)|\ge R$ for any $a$, $b$ and $R$ provided $R$ is large enough. Hence we still have $h(z_n)\to p_{(n_1,\dots,n_{k-1});\pm}$ if $\Im \eta(z_n)\to\pm\infty$,  $z_n\in \cS_{(n_1,\dots,n_{k-1})}$ and $a<\Re\eta(z_n)<b$.

Secondly, for any $n_1,\dots,n_{k-2}\in\Z^{k-2}$ and for any $j_1\in\Z$, the map $  \z_{k-1}=\e(\z_{k-2})$ will map
$$
\{ z \in \cS_{(n_1,\dots,n_{k-2})}|\ \Re\z_{k-2}(z)>\e^{(k-1)}(\alpha), \ |\Im \z_{k-2}(z)-(2j_1+1)\pi|<\pi\}
$$
one-one onto
$$
\tilde\Omega_{k-1}=\C\setminus \left(\{\z_{k-1}\in\R| \z_{k-1}\ge 0\}\cup \{\z_{k-1}|\ |\z_{k-1}|\le \e^{(k-1)}(\alpha)\}\right).
$$
The corresponding curves given by $\Re\eta=\beta$ in $\cS_{(n_1,\ldots,n_{k-2},j_1)}$ and $\cS_{(n_1,\ldots,n_{k-2},j_1-1)}$ give us two branches of curve $\gamma_+$ and $\gamma_-$ satisfying $\Re\eta=\beta$ on $\tilde\Omega_{k-1}\cap \{ \Im\z_{k-1}>0\}$ and $\tilde\Omega_{k-1}\cap \{ \Im\z_{k-1}<0\}$ respectively. Joining the two branches of curve by a compact curve $\gamma$ in $\tilde\Omega_{k-1}$, for instance a circular arc with sufficiently large radius centered at the origin,  gives a subset $U$ with $\partial U=\gamma_+\cup\gamma_-\cup\gamma$ on which $h$ is quasi-conformal. By (4.6), $h$ will maps $\partial U$ to a curve in $\H^2$ such that if $\Im\eta\to\infty$ the image under $h$ will tends to the point $p_{(n_1,\dots,n_{k-2},j_1+1);-}$, and if $\Im\eta\to-\infty$ the image under $h$ will tends to the point $p_{(n_1,\dots,n_{k-2},j_1);+}$. As in the proof of Theorem 3.1, we see that $p_{(n_1,\dots,n_{k-2},j_1+1);-}=p_{(n_1,\dots,n_{k-2},j_1);+} $ which will be denoted by $p_{(n_1,\dots,n_{k-2});j_1}$. It is then not hard to see that
$$
\ol{h(S_{(n_1,\dots,n_{k-1})})}\cap\partial \H^2=\{p_{(n_1,\dots,n_{k-1}) ;j_0} \}_{j_0\in\Z}\cup\{p_{(n_1,\dots,n_{k-1});+},p_{(n_1,\dots,
n_{k-1} );-}\},
$$
and 
$$
\Cal A_0= \{p_{(n_1,\dots,n_{k-2},n_{k-1});j_0} \}_{(n_1,\dots,n_{k-1})\in\Z^{(k-1)},\ j_0\in\Z},
$$ 
is countable and discrete. 
Now for each $(n_1,\dots,n_{k-2})\in\Z^{(k-1)}$ the set 
$$
\{p_{(n_1,\dots,n_{k-2 });j_1}\}_{j_1\in \Z}
$$
is monotone in $j_1$ and we denote $p_{(n_1,\dots,n_{k-2});\pm}=\lim_{j_1\to\pm\infty}p_{(n_1,\dots,n_{k-2});j_1}$.
Since the Hopf differential on $\cS_{(n_1,\dots,n_{k-1})}$ is of the same form (4.3), by the proof of Lemma 3.1, for each $(n_1,\dots, n_{k-1})$, there is a point $z_{(n_1,\dots, n_{k-1})}\in \cS_{(n_1,\dots,n_{k-1})}$ and  two consecutive  points in $\{p_{(n_1,\dots, n_{k-1};j_0)}\}_{j_0\in\Z}$  such that 
$$
\Re \eta(z_{(n_1,\dots, n_{k-1})})=\beta \tag4.7
$$
$$
\lf|\Im  \eta(z_{(n_1,\dots, n_{k-1})})\ri|\le\pi,\tag4.8
$$
and that the distance from $h\lf( z_{(n_1,\dots, n_{k-1})}\ri)$ to the geodesic joining these two consecutive points is bounded by $C_1$ for some  constant $C_1>0$ which is independent of $(n_1,\dots,n_{k-1})$. From (4.7) and (4.8) we have
$$
\lim_{j_1\to\pm\infty}h\lf( z_{(n_1,\dots,n_{k-2}, j_1)}\ri)=p_{(n_1,\dots,n_{(k-2)});\pm}. \tag4.9
$$ 
  Using (4.2) and (4.9), we can argue as before to conclude that
$$
p_{(n_1,\dots,n_{k-3},j_2+1);-}=p_{(n_1,\dots,n_{k-3},j_2) ;+}
$$
which will be denoted by $p_{(n_1,\dots,n_{k-3});j_2}$, and
$$
\ol{h(S_{(n_1,\dots,n_{k-2})})}\cap\partial \H^2=\Cal A_0\cup \left\{p_{(n_1,\dots,n_{k-3},n_{k-2});+}, p_{(n_1,\dots,n_{k-3},n_{k-2} );-}\right\}. 
$$
Let 
$$
\Cal A_1=\{p_{(n_1,\dots,n_{k-2});j_1}\}_{(n_1,\dots,n_{k-2})\in \Z^{k-2},\ j_1\in\Z}.
$$
Then $\Cal A_21$ is countable and each point in $\Cal A_1$ is an isolated accumulation point of $\Cal A_0$. The accumulation points of $\Cal A_1$ are $p_{(n_1,\dots,n_{k-3});j_2}$, $(n_1,\dots,n_{k-3})\in\Z^{(k-2)}$ and $j_2\in\Z$. Continue in this way, we can find $\Cal A_j\subset \partial \H^2$, $0
\le j\le k$ such that each $\Cal A_j$ is countable and discrete for $0\le j\le k-1$ and $\Cal A_j$ consists of all isolated accumulation points of $\Cal A_{j-1}$ for $1\le j\le k$. Moreover,
$$
\ol{h(\C)}\cap\partial \H^2=\cup_{j=0}^k\Cal A_j.
$$
Finally,
We want to prove that $\Cal A_k$ consists of only one point. From the proof, we can see that $\Cal A_k$ consists of at most two points $p$ and $q$ satisfying
$$
\lim_{y\to\infty}h(\i y)=p
$$
and
$$
\lim_{y\to-\infty}h(\i y)=q.
$$
Since
$$
\e^{(k-1)}(t)=\e^{(k-1)}(0)+tg(t)
$$
on $|t|\le 1$, where $g(t)$ is analytic. One can proceed as in the proof of Theorem 4.1 to show that $p=q$ and 
$$
\ol{h(\{z|\ \Re z\le 0 \}}\cap\partial \H^2=\{p\}.
$$
Hence $\Cal A_k$ is a singleton and this completes the proof of the theorem.

\enddemo
\proclaim{Remark 4.2}
The Theorem 4.1 is also true for the Hopf differential 
$$
\e^{(k-1)*}(e^zdz^2)=\e^{(k)}(z)\prod_{j=1}^{k-1}\left[\e^{(j)}(z) \right]^2 dz^2.
$$
In fact, the proof is much easier and can be done by induction since the form of the Hopf differential is not change under the map $\z=e^z$.
\endproclaim
\proclaim{Remark 4.3}
The Theorem 4.1 is not necessary true in general. In fact, it becomes very complicated for the general form as in Theorem 1.2. Even for $\Phi=P(z)\e^{(k)}(z)dz^2$,
the Theorem 4.1 need modification. For instance, if $P(z)=\sqrt{-1}$, then the same argument as in the proof of Theorem 4.1 and using Lemma 4.2 instead of Lemma 4.1 on the region $\{z<\alpha\}$, we see that the set $\Cal A_k$ consists of two points whether than one. So the best to hope for is that $\Cal A_k$ has at most two points for the general form in Theorem 1.2.
\endproclaim

\bigskip
\subheading{\S5 Harmonic diffeomorphisms on hyperbolic plane}

The result in \S2, in particular Proposition 2.1, can be applied to study a conjecture of Schoen, which says that any quasi-symmetric  homeomorphism on $\Bbb S^1$ can be extended to a unique quasi-conformal harmonic diffeomorphism on $\H^2$. The existence part of the conjecture is still open, but there are many partial results, see \cite{Ak, L-T 1--3, T-W 2, S-T-W, H-W, Y}. Schoen's conjecture can be reformulated as follows.  Let $\text {BQD}(\H^2)$ be the space of holomorphic quadratic differentials $\Phi$ on $\H^2$ such that
$$
|||\Phi|||=\sup_{z\in \H^2}||\Phi||(z)<\infty
$$
where $||\Phi||(z)$ is the norm of $\Phi$ at $z$ with respect to the Poincar\'e metric. In \cite{Wn}, the third author proved that for any $\Phi\in\bqd$, there is a unique quasi-conformal harmonic diffeomorphism $u$ on $\H^2$ with $\Phi$ as Hopf differential. This defines a  map $\b$ from $\bqd$ to the universal Teichm\"uller space $\t$ by sending $\Phi$ to the the class of quasi-symmetric homeomorphism containing the boundary value of $u$.  The existence part of the conjecture of Schoen is equivalent to the surjectivity of the map $\b$. Let $\Cal F$ be a subset of $\bqd$, we say that $\b$ is {\it proper} on $\Cal F$ if for any $\Phi_n\in \Cal F$ with $|||\Phi_n|||\to\infty$ we have $d_\t\lf(\b(\Phi_n),0\ri))\to\infty$, where $d_\t$ is the Teichm\"uller metric on $\t$. It is not hard to see that $\b$ is surjective if $\b$ is proper on $\bqd$. It is also not hard to see that if $\b$ is proper on the set of $\Phi\in \bqd$ with $\int_{\H^2}||\Phi||dv_{\H^2}<\infty$   or even on the set $\Phi=\phi dz^2$ with $\phi$ to be a polynomial, then $\b$ is proper on $\bqd$.  Here, we identify $\H^2$ with $\Bbb D$ with the Poincar\'e metric. For the sake of completeness, we give a proof of this fact below. Denote
$$
\Cal F=\left\{\Phi\in\bqd|\ \  \int_{\H^2}||\Phi||dv_{\H^2}<\infty\right\}.
$$
Note that $\Phi=\phi dz^2$, then $\int_{\H^2}||\Phi||dv_{\H^2}=\int_\D|\phi|dxdy$.

\proclaim{Proposition 5.1} 
Let 
$$\Cal G=\{\Phi\in\Cal F|\ \ \Phi=\phi dz^2,\ \text{$\phi$ is a polynomial}\}.$$
Then 
\roster
\item"{(i)}" If $\b$ is proper on $\Cal G$, then $\b$ is proper on $\Cal F$.
\item"{(ii)}" If $\b$ is proper on $\Cal F$, then $\b$ is proper on $\bqd$.
\item"{(iii)}" If $\b$ is proper on $\bqd$, then $\b$ is surjective.
\endroster
In particular, if $\b$ is proper on $\Cal G$, then $\b$ is surjective.
\endproclaim
\demo{Proof} (i) First we prove that if $\b$ is proper on $\Cal G$, then $\b$ is proper on $\Cal F$. Let $\Phi_n\in \Cal F$ such that $|||\Phi_n|||\to\infty$. Suppose that there is a constant $C_1$ such that $d_\t(\b(\Phi_n),0)\le C_1$ for all $n$. Since $\b$ is continuous, there exist $\delta_n>0$ such that if $|||\Phi_n-\Psi|||\le \delta_n$, then
$d_\t(\b(\Psi),0)\le C_1+1$. Hence it is sufficient to prove that $\Cal G$ is dense in $\Cal F$. Let $\Phi=\phi dz^2\in \Cal F$,  then 
$$
\int_{\D}|\phi  |dxdy=\int_{\H^2}||\Phi ||dv_{\H^2}<\infty.
$$
Apply  the mean value inequality on the disk $\D_{z,r}$ with center at $z$ and radius $r=\frac12(1-|z|)$, we can conclude that $||\Phi ||(z)\to 0$ uniformly as $|z|\to 1$. For $0<R<1$, let $\Phi_{R}(z)=\Phi(Rz)$. For any $\epsilon>0$,   we can find $1>\delta>0$ such that if
$1-\delta\le |z|<1$ then $||\Phi||(z)\le \frac12\epsilon$.   Then for $1-\frac12\delta\le |z|<1$ and for
$R$ large enough, so that  $R|z|\ge 1-\delta$
$$
||\Phi_R(z)||=\frac{(1-| z|^2)^2}{(1-|Rz|^2)^2}||\Phi(Rz)||\le \frac 12\epsilon.
$$
On the other hand, for $|z|\le 1-\frac12\delta$, $\phi_R(z)\to\phi(z)$ uniformly,
as $R\to 1$. Hence we can find $R$ large enough, so that
$$
||\Phi_R(z)-\Phi(z)||\le   \epsilon
$$
for all $z\in\D$. Hence $|||\Phi_R -\Phi |||_{BQD}\le \epsilon.$ But $\Phi_R$ is analytic
on $|z|<\frac1R$ which is large than 1. So it can be approximated uniformly on $\D$ by polynomials. This completes the proof of (i)

(ii) We will prove that if $\b$ is proper on $\Cal F$, then $\b$ is proper on $\bqd$. Let $\Phi\in \Cal F$ and let $\b(\Phi)=[f]$ where $f$ is a quasi-symmetric homeomorphism of $\Bbb S^1$ fixing $1,\ i, \ -i$. Let $d_\t([f],0)=C_1$. Then there exist smooth quasi-symmetric functions $g_k$ fixing $1,\ i,\ -i$ such that $g_k\to f$ in $C^\alpha$ norm for some $1>\alpha>0$ and such that $d_\t([g_k],0)\le C_2$ which depends only on $C_1$. Moreover,  $C_1\to\infty$ if and only if $C_2\to\infty$. These follow from theorem 2 and remark (1) in  \cite{D-E}. By theorem 6.4 in \cite{L-T 3}, see also \cite{T-W 2}, for each $k$ there exists a unique $\Psi_k\in\bqd$ such that $\b(\Psi_k)=[g_k]$ with $\Psi_k\in\Cal F$. By the assumption, we have $|||\Psi_k|||\le C_3$ for all $k$, where $C_3$ depends only on $C_1$. Note that $\Psi_k$ is the Hopf differential of quasi-conformal harmonic diffeomorphism on $\H^2$ with boundary value $g_k$. Hence $\Psi_k(z)\to\Phi(z)$ for all $z\in \D$ and so
$$
\limsup_{k\to\infty}|||\Psi_k|||\ge |||\Phi|||.
$$
From this, it is easy to see that $\b$ is proper on $\bqd$.

(iii) We will prove that if $\b$ is proper on $\bqd$, then $\b$ is surjective. Let $[f]$ be a class of quasi-symmetric homeomorphism on $\Bbb S^1$ such that $[f]$ is in the closure of $\b(\bqd)$. Then there exists   $f_n$  quasi-symmetric homeomorphisms on $\Bbb S^1$ fixing $1,\ i,\ -i$ such that $f_n\to f$ uniformly, and $[f_n]=\b(\Phi_n)$. Since $[f_n]$ are uniformly bounded on $\t$, $\Phi_n$ are uniformly bounded in $\bqd$. By theorem 13 in \cite{Wn}, the quasi-conformal harmonic diffeomorphisms $u_n$ with Hopf differentials $\Phi_n$ has complex dilatation $\mu_n$ satisfying $|\mu_n|\le \mu<1$ for some constant $\mu$ independent of $n$. Passing to a subsequence if necessary,   $u_n$ converges uniformly on $\overline\D$ to a quasi-conformal harmonic diffeomorphism on $\H^2$ with boundary value $f$. Hence $[f]$ is in $\b(\bqd)$. Combine with the theorem 4.1 in \cite{T-W 2}, we conclude that $\b$ is surjective.
   \enddemo

\proclaim{Proposition 5.2} Let $\Phi_n\in \bqd$ satisfying $\int_{\H^2}||\Phi_n||<\infty$ and $|||\Phi_n|||\to\infty$. Suppose for all $k>0$, 
$$
\lim_{n\to\infty}\frac{\int_{U_n}||\Phi_n||dv_{\H^2}}{\int_{\H^2}||\Phi_n||dv_{\H^2}}=0
$$
where $$
U_n=\{z\in \H^2|\ ||\Phi_n||(z)\le k|||\Phi_n|||^\frac34\}.
$$
Then $d_{\t}(\b(\Phi_n),0)\to\infty.$
\endproclaim

\demo{Proof} Identify $\H^2$ with the unit disk $\D$ equipped with the Poincar\'e metric. Let $\Phi_n= \phi_n dz^2$. Let $\mu_{n,0}$ be the supremum of the modulus of the  complex dilation of $\b(\Phi_n)$, then
$$
d_{\t}(\b(\Phi_n),0)=\frac12\log\frac{1+\mu_{n,0}}{1-\mu_{n,0}}.
$$ 
Since 
$$
\int_{\D}|\phi_n|dxdy=\int_{\H^2}||\Phi_n||dv_n<\infty,
$$
we can apply   the main inequality in \cite{R-S} and conclude 
that 
$$
\frac{1}{1+\mu_{n,0}}\int_{\H^2}||\Phi_n||dv_{\H^2}\le
\int_{\H^2}\frac{1}{1+|\mu_n|}||\Phi_n|| dv_{\H^2} \tag5.1
$$
where $\mu_n$ is the complex dilatation of the  quasi-conformal harmonic diffeomorphism with Hopf differential $\Phi_n$. Let $1>\delta>0$ and $k>1$ be fixed numbers. Define
$$
 D_n=\{z\in \H^2|\ ||\Phi_n||(z)\ge k||\Phi_n||_{BQD}^\frac34\},
$$
$$
U_n=\{z\in \H^2|\ ||\Phi_n||(z)< k||\Phi_n||_{BQD}^\frac34\}.
$$

We have 
$$
\int_{\H^2} \frac{1}{1+|\mu_n|}||\Phi_n||dv_{\H^2}
 =\lf(\int_{D_n} +\int_{U_n}\ri)\frac{1}{1+|\mu_n|}||\Phi_n||
dv_{\H^2}.\tag5.2
$$
By Proposition 2.1, for any $z\in D_n$, the maximal $\Phi_n$-radius $R_{n,z}$ satisfies $R_{n,z}\ge C_1k$ for some absolute constant $C_1>0$. By the result on page 63 of \cite{Hn}, we have $|\mu_n|\ge \eta(k)$ for some constant $\eta(k)$ such that $\eta(k)\to1$ as $k\to\infty$. Hence  we have
$$
\int_{D_n} \frac{1}{1+|\mu_n|}||\Phi_n||dv_{\H^2}\le
\frac{1}{1+\eta(k)}\int_{D_n}||\Phi_n||dv_{\H^2}.\tag5.3
$$
For any $0<\delta<1$,  there exists $n_0$ such that if $n\ge n_0$ then 
$$
\int_{U_n}||\Phi_n||dv_{\H^2}\le \delta\int_{\H^2}||\Phi_n||dv_{\H^2}.
$$
Combine this with (5.1)--(5.3), we have for $n\ge n_0$,

$$
\frac{1}{1+\mu_{n,0}}\int_{\H^2}||\Phi_n||dv_{\H^2}\le
\frac{1}{1+\eta(k)}\int_{D_n}||\Phi_n||dv_{\H^2}
+\delta\int_{\H^2}||\Phi_n||dv_{\H^2}.
$$
Hence
$$
\frac{1}{1+\mu_{n,0}}\le \frac{1}{1+\eta(k)}+\delta.
$$
Let $k\to\infty$, and then let $\delta\to0$, we have
$$
\limsup_{n\to\infty}\frac{1}{1+\mu_{n,0}}\le\frac12.
$$
Since $0\le \mu_{n,0}<1$ for all $n$, we have $\lim_{n\to\infty}\mu_{n,0}=1$. From this, it is easy to see that $d_{\t}(\b(\Phi_n),0)\to\infty.$
\enddemo

\proclaim{Corollary 5.1}  Let $\Cal F$ as in Proposition 5.1.  Let $\Phi_n=\phi_ndz^2\in \Cal F$. Suppose
$|||\Phi_n|||=1$. Let $A_n=\int_{\H^2}||\Phi_n||dv_{\H^2}=\int_\D |\phi_n|dxdy$. Suppose
$\lim_{n\to\infty}\phi_n/A_n=\psi$  such that
$\int_{\D}|\psi|dxdy=1$, then  $\b(t_n\Phi_n)\to\infty$  for any $t_n\to\infty$. In particular $\b$ is proper on any finite dimensional subspace of    $\Cal F$.
 \endproclaim

\demo{Proof}  For any $1>\delta>0$, we can find $1>r_0>0$, such that
$$
\int_{\D_{r_0}}|\psi|dxdy=1-\delta, 
$$
where $\D_{r_0}=\{|z|<r_0\}$. We have
$$
\int_{\D_{r_0}}\frac{|\phi_n|}{ A_n}\ge 1-2\delta,\tag5.4
$$
provided $n$ is large enough. For any $k>0$, let $\epsilon_n=kt_n^{-\frac14}$, where $t_n\to\infty$.  Let
$$
U_n=\{z\in\H^2|\
||t_n\Phi_n||(z)\le k|||t_n\Phi_n|||^\frac34 \}=\{z\in\H^2|\
||\Phi_n||(z)\le \epsilon_n \}.
$$
Note that if $\Phi_n=\phi_n dz^2$, then for $z\in U_n\cap\D_{r_0}$
$$
|\phi_n|(z)\le C_1(1-r_0)^{-2}||\Phi_n||(z)\le \epsilon_n (1-r_0)^{-2}
$$
for some absolute constant $C_1$. Hence by (5.4), 
$$
\split
\int_{U_n}||\Phi_n||dv_{\H^2}&\le \int_{U_n\cap\D_{r_0}}|\phi_n|dxdy+\int_{\D\setminus\D_{r_0}} |\phi_n|dxdy\\
&\le  \frac{\pi C_1\epsilon_n}{(1-r_0)^2}+2\delta A_n .
\endsplit\tag5.5
$$
Since $|||\Phi_n|||=1$, by applying the mean value inequality to $|\phi_n|(z)$  at a point $z_n$ with $||\Phi_n||(z_n)\ge\frac12$, we conclude that $A_n\ge C_2$ for some abolute constant $C_2>0$.  By the definition of $\epsilon_n$, we have $\epsilon_n\to0$ as $n\to\infty$. Hence (5.5) implies that we
have
$$
\split
\limsup_{n\to\infty}\frac{\int_{U_n}||t_n\Phi_n||dv_{\H^2}}{\int_{\H^2}||t_n\Phi_n|dv_{\H^2}}
&=\limsup_{n\to\infty}\frac{\int_{U_n}||\Phi_n||dv_{\H^2}}{A_n}\\ 
&\le 2\delta.
\endsplit
$$
Since $\delta$ is arbitrary, $t_n\Phi_n$ satisfies the conditions in the Proposition 5.2. Hence the first part of corollary is proved.

To prove the second part of the corollary, let $\Cal H$ be a finite dimensional subspace of $\Cal F$ with basis $\Psi_1,\dots,\Psi_k$. If $\Phi_n\in\Cal H$ is  such that $|||\Phi_n|||\to\infty$, then $\Phi_n=t_n\tilde\Phi_n$ such that $|||\tilde\Phi_n|||=1$ and $\tilde\Phi_n=\sum_{j=1}^ka_{n,j}\Psi_j$ with $\sum_{j=1}^ka_{n,j}^2$ being uniformly bounded. From this, it is easy to see that the second part of the corollary follows.
\enddemo
 Let $\Phi=\phi dz^2\in \Cal F$. Define
$$
|\Phi|_\infty=\inf_{a\in \D}\sup_{z\in \D}|\phi(z)|\frac{|1-\bar a
z|^4}{(|1-|a|^2)^2}=\inf_{a\in \D}\sup_{\zeta\in \D}|\tilde
\phi_a(\zeta)|  ,
$$
where $\tilde \phi_a d\zeta^2=\Phi$ and $\zeta=\frac{z-a}{1-\bar az}$.

$$
|\Phi|_{L^1}=\int_\D|\phi| dxdy=\int_{\H^2}||\Phi||dv_{H^2}.
$$
Note that  
$$
\int_\D|\phi|dxdy=\int_\D|\tilde\phi_a| dxdy
$$
if $\phi$ and $\tilde \phi_a$ are related as above.

\proclaim{Corollary 5.2} Let $\Cal F_1=\{\Phi\in BQD|\ |||\Phi|||=1,
|\Phi|_\infty<\infty\}$. Let $\Phi_n\in \Cal F_1$, and $t_n\to\infty$ be a sequence such that $ |\Phi_n|_\infty   |\Phi_n|^{-2}_{L^1}=o(t^{\frac14}_n)$, then    $d_{\Cal T}(\b(t_n\Phi_n),0)\to\infty$. In particular, if
$|\Phi_n|_\infty   |\Phi_n|^{-2}_{L^1}\le C$ for some constant $C$ independent of $n$, then $d_{\Cal T}(\b(t_n\Phi_n),0)\to\infty$, for any $t_n\to\infty$.
\endproclaim
\demo{Proof} For each $n$, by the definition, by a linear fractional
transformation of $\D$ if necessary, we may assume that $\Phi_n=\phi_n dz^2$,
with $\lf(\sup_{z\in \D}|\phi_n|(z)\ri) |\phi_n|^{-2}_{L^1}=o(t_n^\frac14)$. Let
$M_n=\sup_{z\in
\D}|\phi_n|(z)$ and $I_n=|\phi_n|_{L^1}$. We claim that 
$$
I_n^2\le C_1M_n\tag5.6
$$
for all $n$ and for some absolute constant $C_1$.  Fix   $n$, take $r_0$ such that
$$
\int_{\D(r_0)}|\phi_n|dxdy=\frac12I_n,
$$
where $\D(r_0)=\{z|\ |z|<r_0\}$. Since $|||\Phi_n|||=1$, we have
$$
\frac12I_n=\int_{\D(r_0)}|\phi_n|dxdy\le C_2\int_0^{r_0}(1-r)^{-2}\le \frac{C_2}{1-r_0}\tag5.7
$$
where $C_2$ is an absolute constant. On the other hand
$$
\frac12I_n=\int_{\D\setminus\D(r_0)}|\phi_n|dxdy\le M_n\cdot2\pi(1-r_0).
$$
Hence
$$
\split
\frac12I_n&\le \frac{C_2}{1-r_0}\\
&\le \frac{2\pi M_nC_2}{\frac12I_n}.
\endsplit
$$
From this (5.6) follows.

Now for any $k>0$, let 
$$
\split
U_n&=\{z\in \H^2|\ ||t_n\Phi_n||(z)\le k|||t_n\Phi_n|||^{3/4}\}\\
&= \{z\in \H^2|\ ||\Phi_n||(z)\le k\epsilon_n|||\Phi_n|||^{3/4}\}
\endsplit
$$
where $\epsilon_n=t_n^{-\frac14}$.  Let $\delta_n=\lf(\frac{k\epsilon_n}{M_nI_n^{-2}}\ri)^\frac12.$ By (5.6) and the fact that $\epsilon_n\to0$, we have $\delta_n\to0$ as $n\to\infty$. As in the proof of Corollary 5.1, $I_n\ge C_3$ for some absolute constant $C_3$. Hence $\delta_nI_n^{-1}\to0$ as $n\to\infty$. If $n$ is large enough so that $\delta_nI_n^{-1}<1$, then we take $r_n$ such that $1-r_n=
\delta  I_n^{-1}$. 
Denote
$\D(r_n)=\{ z|\ |z|<r_n\}$. Then
$$
\split
\int_{(\D\setminus\D_n)\cap U_n}|\phi_n|dxdy&\le 2\pi M_n(1-r^2_n)\\
&\le 4\pi(M_nI_n^{-2})\delta_n I_n.
\endsplit\tag5.8
$$
On the other hand, as in the proof of (5.7), we have 
$$
\split
\int_{\D(r_n)\cap U_n}|\phi_n|dxdy&\le C_4k\epsilon_n(1-r_n)^{-1}\\
&=C_4k\epsilon_n\delta^{-1}_nI_n 
\endsplit
$$
for some constant $C_4$ independent of $n$.
Hence
$$
\split
\int_{U_n}|\phi_n|dxdy&\le
C_5\lf\{M_nI_n^{-2}\delta_n+k\epsilon_n\delta_n^{-1}\ri\}I_n\\
&=2C_5\lf\{ k\epsilon_nM_nI_n^{-2}\ri\}^\frac12I_n 
\endsplit
$$
for some constant $C_5$ independent of $n$. 
Hence if $n$ is large enough, 
$$
\frac{\int_{U_n}|\phi_n|dxdy}{\int_{\D}|\phi_n|dxdy}\le 2C_5\lf\{ k\epsilon_nM_nI_n^{-2}\ri\}^\frac12.
$$
By the assumption, the right side of the above inequality tends to zero as $n\to0$ and the corollary follows.
\enddemo
  
{\bf Example 1.} Let $\Phi_n =\phi_n dz^2=c_nn^2 z^ndz^2$, where $c_n$ is chosen so that $|||\Phi_n|||=1$. Direct computations show  that   $C^{-1}\le c_n\le C$ for some positive constant $C$  independent of $n$.
Then $\sup_{z\in \D}|\phi_n|(z)=c_n n^2$, and $\int_{\D}|\phi_n|dxdy=\frac{2\pi n^2
c_n}{n+2}$. Hence the $|\Phi_n|_\infty |\Phi_n|^{-2}_{L_1}\le C'$ for some constant $C'$ independent of $n$. Hence by Corollary 5.2, for any subsequence $n_k$ and for any $t_k\to\infty$ we have $d_{\Cal T}(\b(t_k\Phi_{n_k}),0)\to\infty$.\vskip .2cm

 {\bf Example 2.}  Consider
$(z-1)^ndz^2$. Then direct computation shows
$$
 \sup_{z\in \D}(1-|z| )^2|z-1|^n =c_nn^{-2}2^n 
 $$
where $C^{-1}\le c_n\le C$ for some constant $C>0$ independent of $n$. Let
$\Phi_n=c_n'n^22^{-n}(z-1)^ndz^2=\phi_n(z)dz^2$, where $c_n'$ is chosen so
that $|||\Phi_n|||=1$. Note that $\tilde C^{-1}\le c_n'\le \tilde C$ for some $\tilde C>0$
independent of $n$. Let $a=-\frac{n+8-4\sqrt{n+4}}{n}$, note that 
$-1<a<0$. 
$$
\split
\sup_{z\in \D}|\phi_n(z)|\frac{|1-az|^4}{(1-a^2)^2}
&=c_n'n^22^{-n}\sup_{0\le \theta\le
2\pi}|e^{i\theta}-1|^n\ |1-ae^{i\theta}|^4(1-a^2)^{-2}\\
&=c_n'n^2 \sup_{0\le \theta\le
2\pi}\sin^n\frac\theta2\lf [(1-a)^2+4a\sin^2\frac\theta 2\ri]^2(1-a^2)^{-2}\\
&=c_n'n^2\sup_{0\le t\le 1}t^n\lf [(1-a)^2+4at^2\ri]^2(1-a^2)^{-2}.
\endsplit\tag5.9$$
Let $f(t)=t^n\lf [(1-a)^2+4at^2\ri]^2$, then $f(0)=0$. Suppose $f(t)$, $0\le
t\le 1$ attains its maximum at $t_0\in (0,1)$, then $f'(t_0)=0$, and
$$
nt_0^{n-1} \lf [(1-a)^2+4at_0^2\ri]^2+t_0^n\cdot 2\lf
[(1-a)^2+4at_0^2\ri] \cdot 8at_0=0.$$
Hence
$t_0^2=-\frac{n(1-a^2)}{4a(n+4)}=1$ by the choice of $a$, which is
impossible. Hence, for $0\le t\le 1$, $f(t)$ attains its maximum at $t=1$. By (5.9), we have
$$
\split
|\Phi_n|_\infty&\le c_n'n^2\lf[(1-a)^2+4a\ri]^2(1-a^2)^{-2}\\
&=c_n'n^2(1+a)^2(1-a)^{-2}.
\endsplit
$$
Since $a<0$, $1-a>1$. Also 
$$
1+a=1-\frac{n+8-4\sqrt{n+4}}{n}=\frac{-8+4\sqrt{n+4}}{n}.
$$
Hence 
$$
|\Phi_n|_\infty\le C_1n \tag5.10$$
for some constant $C_1$ independent of $n$.
On the other hand
$$
\split
\int_\D|z-1|^ndxdy&=\int_{2\pi}^{0}\int_0^1|r^2+1-2r\cos\theta|^{\frac
n2}rdrd\theta\\
&=\int_{-\pi}^{\pi}\int_0^1|r^2+1+2r\cos\theta|^{\frac
n2}rdrd\theta\\ &\ge\int_0^{\frac{1}{\sqrt
n}}\int_0^1(r+1)^n|1-\frac{r\theta^2}{(1+r)^2}|^\frac n2rdrd\theta\\
&\ge  \frac{C_2}{\sqrt n}\int_0^1 (r+1)^n\lf (1-\frac
n2\frac{r}{n(1+r)^2}\ri)rdr\\
&\ge C_32^nn^{-\frac32}\\
&\ge C_4n^{\frac12} 
\endsplit
$$
for some constants $C_2--C_4$ independent of $n$.
By Corollary 5.2, we also have $\b(t_k\Phi_{n_k})\to\infty$ for all subsequence $n_k$ and for all $t_k\to\infty$.\vskip .2cm

In the last section of \cite{Wn}, it was proved that $\b$ is continuous, and in section 4 of \cite{T-W 2}, it was proved that the image of $\b$ is open and $\b$ is a diffeomorphism from $\bqd$ into $\t$. From the proof of the proposition 14 in \cite{Wn},  $\b$ is in fact uniformly continuous on bounded subsets of $\bqd$. On the other hand, we have the following:
\proclaim{Proposition 5.3} Let $R>0$ and let 
$$
B(R)=\{\Phi\in\bqd|\ |||\Phi|||\le R\}.
$$
Let $\t(R)=\b\lf(B(R)\ri)$. Then $\b^{-1}$ is uniformly continuous on $\t(R)$. 
\endproclaim
\demo{Proof} For any complex measurable function $\mu$ on $\D$ such that $||\mu||_\infty<1$, denote $F^\mu$ to be the unique quasi-conformal map on $\D$ with boundary value $f^\mu$ which fixes $1,\ i, \ -i$. Suppose $f^\mu$ can be extended to a quasi-conformal harmonic diffeomorphism, then the harmonic map will be denoted by $\tf^\mu$ and its complex dilation is denoted by $\tm$. By theorem 13 in \cite{Wn}, there exists $0<k<1$ such that if $[f]\in \t(R)$ then $||\mu||_\infty\le k$, where $\mu$ is the complex dilatation of an extremal quasi-confomal map with boundary value in $[f]$.  
$$
\t^*(R)=\{\mu|\ \text{$\mu$ is measurable, $||\mu||_\infty\le k$ and $[f^\mu]\in \t(R)$}\}.
$$
Note that if $\mu\in \t^*(R)$, then $f^\mu$ can be extended to a quasi-conformal harmonic diffeomorphism with Hopf differential in $B(R)$. We claim that for any $\epsilon>0$, there is $\delta>0$ such that if $\mu$, $\nu$ in $\t^*(R)$ and $||\mu-\nu||_\infty\le \delta$ then $||\tm-\tn||_\infty\le\epsilon$. If the claim is true, then by the definition of $d_\t$ and by the method 
as in the proof of proposition 14 in \cite{Wn}, one can conclude that $\b^{-1}$ is uniformly continuous on $\t(R)$.

First we prove the following, given $\epsilon>0$, there is $\delta>0$ such that if $\mu$ and $\nu$ are in $\t^*(R)$, then $|\tm(0)-\tn(0)|\le \epsilon$. Suppose not, then there  is  $\epsilon>0$  and two sequences $\mu_n$, $\nu_n$ in $\t^*(R)$ such that $||\mu_n-\nu_n||_\infty\to0$, but $|\tm_n(0)-\tn_n(0)|\ge \epsilon$. Since $||\tm_n||_\infty\le k_1$ and $||\tn_n||_\infty\le k_1$ for some $0<k_1<1$ by \cite{Wn},  passing to   subsequences if necessary, $\tf^{\mu_n}$ and $\tf^{\nu_n}$ converge uniformly on $\overline{\D}$ to normalized  quasi-conformal harmonic diffeomorphisms $H_1$ and $H_2$ respectively. Since $||\mu_n-\nu_n||_\infty\to0$ and $\mu_n$, $\nu_n$ are in $\t^*(R)$, $H_1$ and $H_2$ must have the same boundary value and so $H_1=H_2$ by \cite{L-T 3}. It then follows that $|\tm_n(0)-\tn_n(0)|\to0$, which is a contradiction. 

Now, for any $\epsilon>0$, let $\delta>0$ be as above. Let $\mu$, $\nu$ be in $\t^*(R)$ such that $||\mu-\nu||_\infty\le\delta$. Let $a\in \D$ and $\phi (z)=(z-a)/(1-\bar az)$. Define $\mu_1$ and $\nu_1$ by
$$
\mu_1(\phi(z))=\mu(z)\lf(\frac{\phi'(z)}{|\phi'(z)|}\ri)^2
$$ and 
$$\nu_1(\phi(z))=\nu(z)\lf(\frac{\phi'(z)}{|\phi'(z)|}\ri)^2.
$$ 
Then $f^{\mu_1}=h^\mu\circ f^\mu\circ\phi^{-1}$, and $f^{\nu_1}= h^\nu\circ f^\nu\circ\phi^{-1}$ where $h^\mu$ and $h^\nu$ are the linear fractional transformations which map $\D$ onto itself  and are chosen so that   $f^{\mu_1}$ and $f^{\nu_1}$ fix $1,\ i, \ -i$ respectively.
Obviously $f^{\mu_1}$ and $f^{\nu_1}$ have quasi-conformal representatives $\tf^{\mu_1}$ and $\tf^{\nu_1}$. In fact,
$$
\tf^{\mu_1}=h^\mu\circ \tf^\mu\circ\phi^{-1}
$$
and
$$
\tf^{\nu_1}=h^\mu\circ \tf^\nu\circ\phi^{-1}.
$$
Moreover, the Hopf differentials of $\tf^{\mu_1}$ and $\tf^{\nu_1}$ are in $B(R)$. Hence $\mu_1$, $\nu_1$ are in $\t^*(R)$ becasue $||\tm_1||_\infty=||\mu||_\infty\le k$ and $||\tn_1||_\infty=||\nu||_\infty\le k$ .

We also have
$$
\tm_1(\phi(z))=\tm(z)\lf(\frac{\phi'(z)}{|\phi'(z)|}\ri)^2
$$ 
and 
$$\tn_1(\phi(z))=\tn(z)\lf(\frac{\phi'(z)}{|\phi'(z)|}\ri)^2.
$$ 
 Note that  
$$
||\mu_1-\nu_1||_\infty=||\mu-\nu||_\infty\le\delta.
$$
Therefore
$$
|\mu(a)-\nu(a)|=|\mu_1(0)-\nu_1(0)|\le \epsilon.
$$
Since $a$ is any point in $\D$, the claim follows.
\enddemo

\bigskip
%
%
%
%
%
%
%
\input epsf.tex  
\def\caption#1{\hfill \\ vskip -3ex plus 1ex minus 0ex
  \hbox{}\hfil{\footnotesize #1}\hfil\hbox{}}
\def\ifnextchar#1#2#3{\let\tmpnce=#1%
    \def\tmpnca{#2}\def\tmpncb{#3}\futurelet\tmpncc\ifnch}%
\def\ifnch{\ifx\tmpncc\tmpnce\let\tmpncd=\tmpnca%
    \else\let\tmpncd=\tmpncb\fi\tmpncd}
%
\def\tpeinture #1 by #2 (#3){
  \vtop to #2{
    \hrule width #1 height 0pt depth 0pt
    \vfill
    \epsfysize=#2 \epsfbox{#3}
    }
  }
\def\xpeinture #1 by #2 (#3){
  \hbox{$\vcenter to #2{
    \hrule width #1 height 0pt depth 0pt
    \vfill
    \epsfxsize=#1 \epsfbox{#3}
    }$}
  }
\def\ypeinture #1 by #2 (#3){
  \hbox{$\vcenter to #2{
    \hrule width #1 height 0pt depth 0pt
    \vfill
    \epsfysize=#2 \epsfbox{#3}
    }$}
  }
\def\peinture{\ypeinture}
\def\bpeinture #1 by #2 (#3){
  \vbox to #2{
    \hrule width #1 height 0pt depth 0pt
    \vfill
    \epsfysize=#2 \epsfbox{#3}
    }
  }
\def\scaledpik[#1] #2 by #3 (#4){{ %
   \dimen0=#2 \dimen1=#3 %
   \if#1t
       \tpeinture \dimen0 by \dimen1 (#4)%
   \else\if#1b
        \bpeinture \dimen0 by \dimen1 (#4)%
        \else\if#1x
             \xpeinture \dimen0 by \dimen1 (#4)%
             \else
             \ypeinture \dimen0 by \dimen1 (#4)%
            \fi\fi
    \fi}}
\def\scaledpic #1 by #2 (#3){{%
   \dimen0=#1 \dimen1=#2%
   \peinture \dimen0 by \dimen1 (#3) %
   }}
\def\scaledpicture{\ifnextchar[{\scaledpik}{\scaledpic}}
\def\centredpicture #1 by #2 (#3){
   \par\centerline{\hbox{
   \scaledpicture #1 by #2 (#3)}
   }}

\subheading{Appendix: Trajectories and Image Accumulation}

In this appendix, seven figures of horizontal trajectories defined
by a quadratic differential are shown.  These pictures of
trajectories are produced by programming in Mathematica.  Some
trajectories may be broken due to slow convergence of the algorithm.
In fact it should be smoothly defined for all time.  Nevertheless,
the qualitative behavior of the trajectory patterns is shown
clearly.
In some figures, the correponding harmonic map produces an image
which has a good accumulation structure on the
boundary.  This structure is also shown on the unit disk.

\head Finitely Many Accumulations
\endhead

\subsubhead Figure 1. $\Phi=e^z dz^2$
(See Theorem~3.1 and Corollary~3.1)
\endsubsubhead
This example is the basis of all others.  The trajectories have a $2\pi
i$ periodicity.
The image of the corresponding harmonic map has an accumulation point at $-1$.

\centredpicture 11cm by 7cm (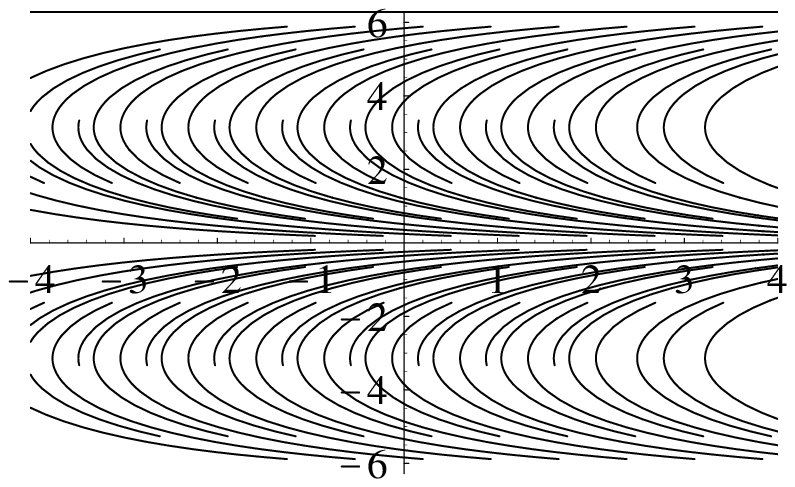)
\par\vskip 3ex plus 1ex minus 1ex\par
\centredpicture 7cm by 7cm (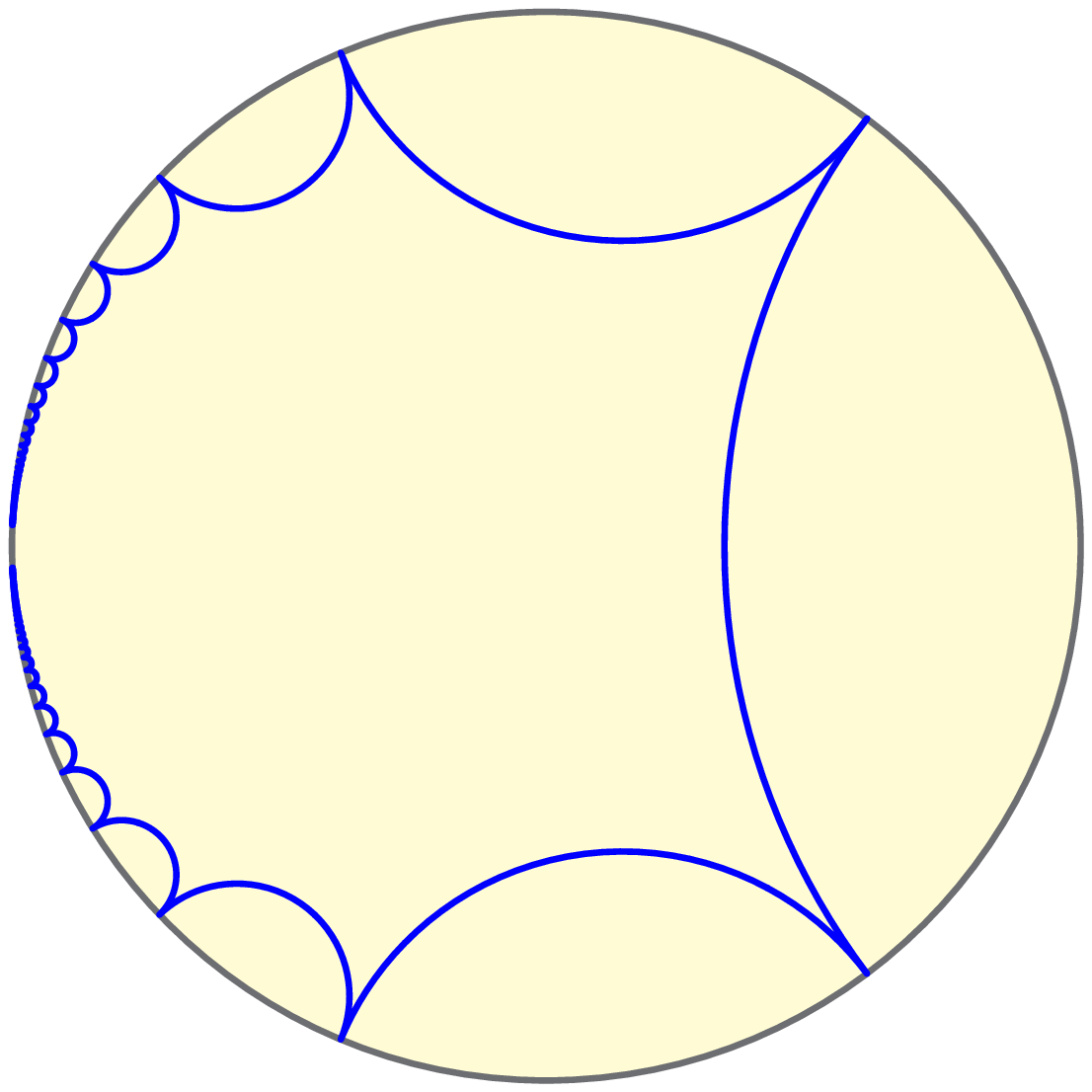)

\filbreak

\subsubhead Figure 2. $\Phi=(e^z+1) dz^2$
(See Theorem~4.1)
\endsubsubhead
\centredpicture 11cm by 7cm (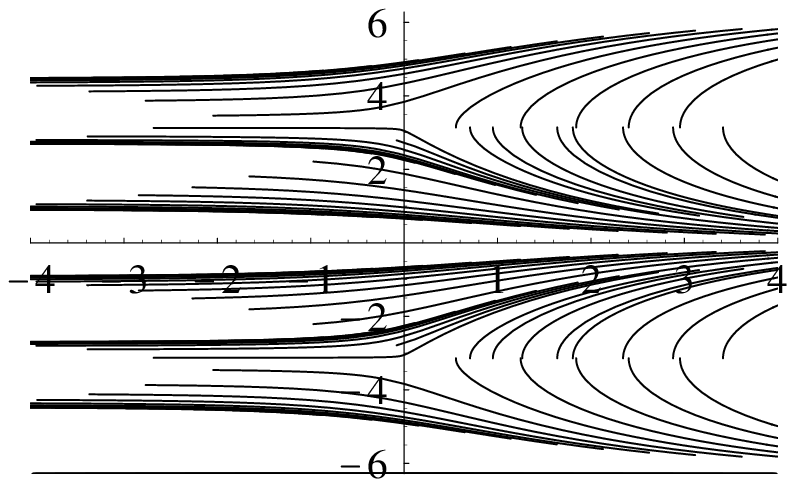)

\break
\subsubhead Figure 3. $\Phi=(e^z-1) dz^2$
(See Theorem~4.1)
\endsubsubhead
From this and the previous one, a lower order term may
significantly change the behavior of the harmonic map.
This one has two accumulation points at $\pm 1$ while the
previous one has only one.

\centredpicture 11cm by 6.5cm (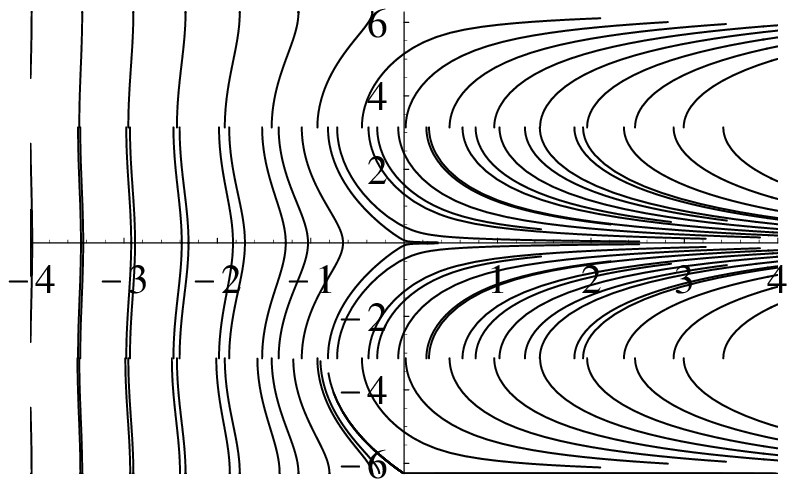) 
\par\vskip 3ex plus 1ex minus 1ex\par
\centredpicture 6.5cm by 6.5cm (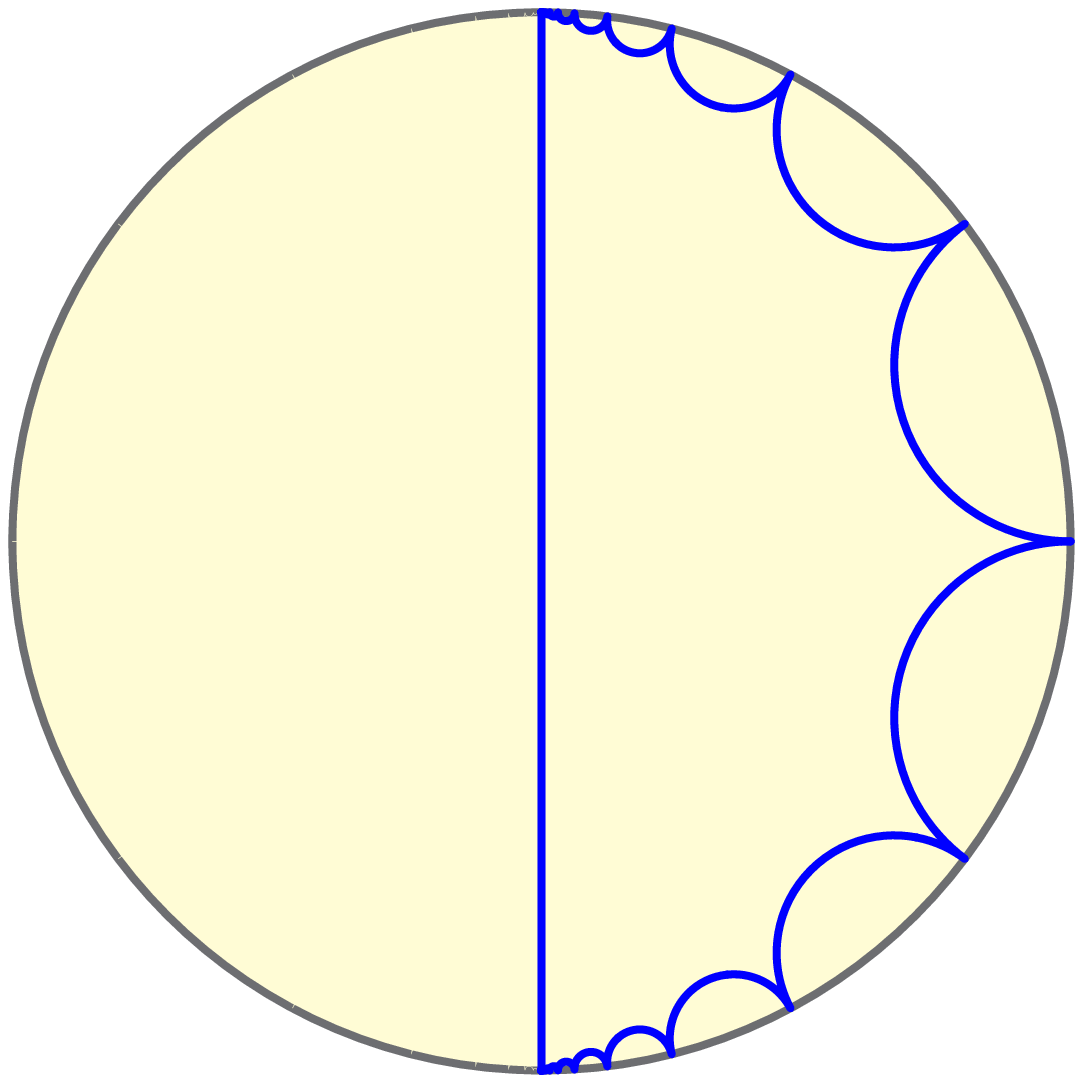)

\filbreak

\subsubhead Figure 4. $\Phi=\sinh^2 z dz^2$
(See Theorem~4.1)
\endsubsubhead
This is another example that the fundamental region is
different while there are also two accumulation points.

\centredpicture 11cm by 6.5cm (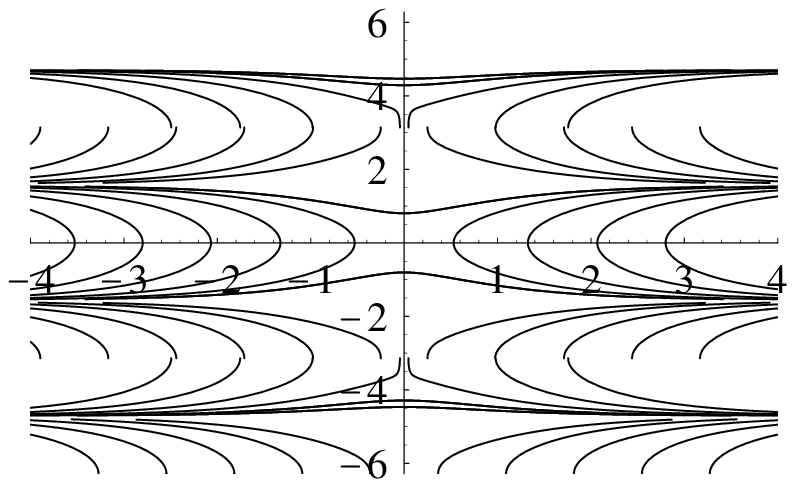) 
\par\vskip 3ex plus 1ex minus 1ex\par
\centredpicture 6.5cm by 6.5cm (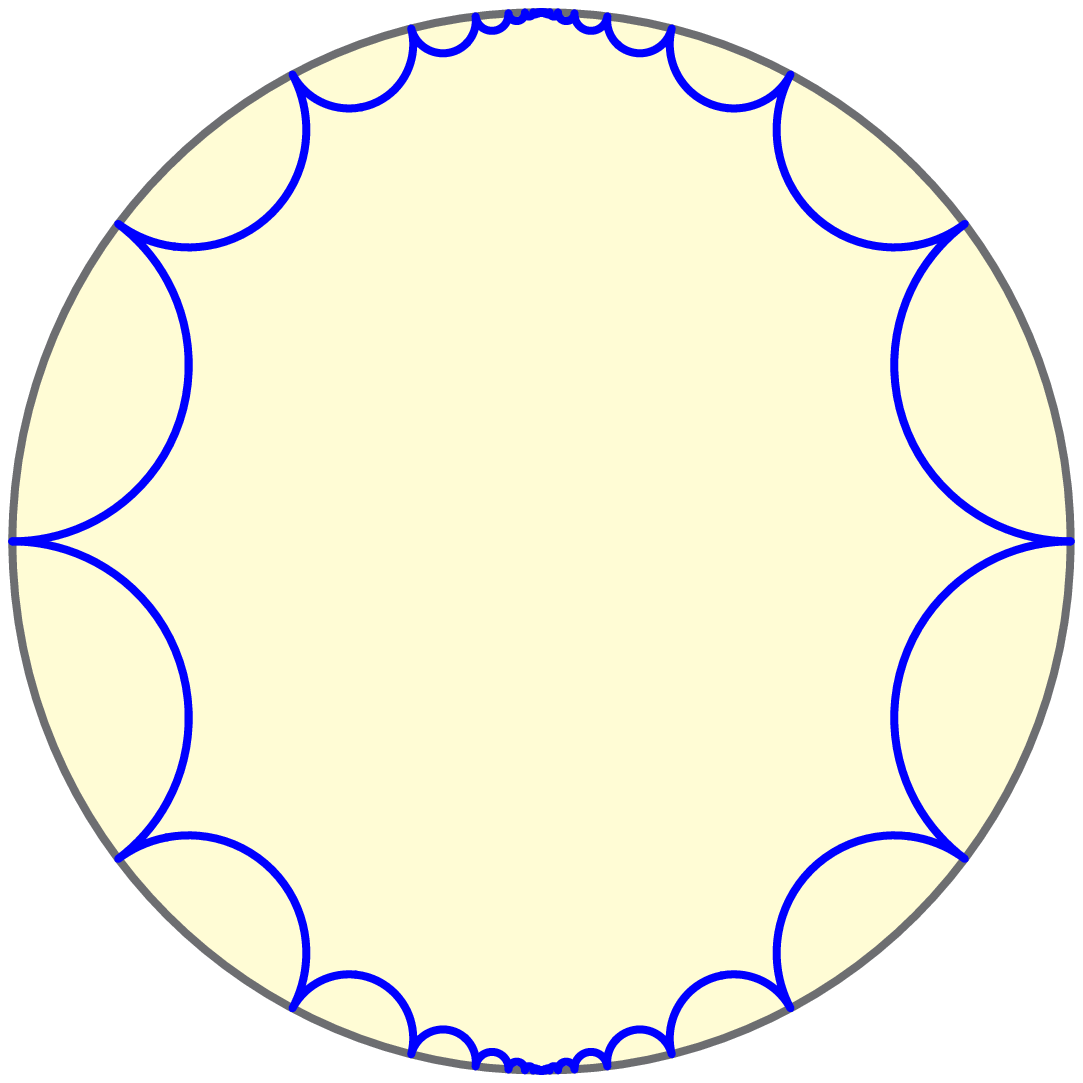)

\filbreak

\subsubhead Figure 5. $\Phi=(z^2-1)e^{z^2} dz^2$
(See Theorem~3.1 and Corollary~3.1)
\endsubsubhead
The image of this example should have 2 accumulation
points.
\centredpicture 11cm by 7cm (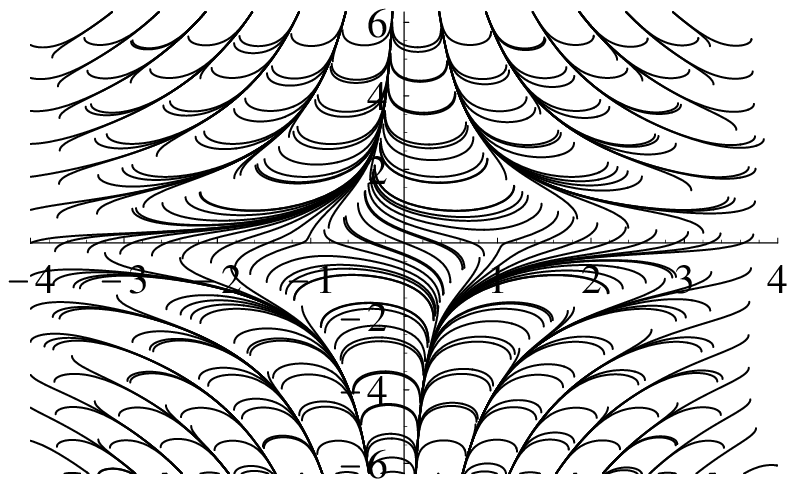)

\subsubhead Figure 6. $\Phi=e^{z^3+z^2} dz^2$
(See Theorem~3.1 and Corollary~3.1)
\endsubsubhead
The image of this example should have 3 accumulation
points.
\centredpicture 11cm by 7cm (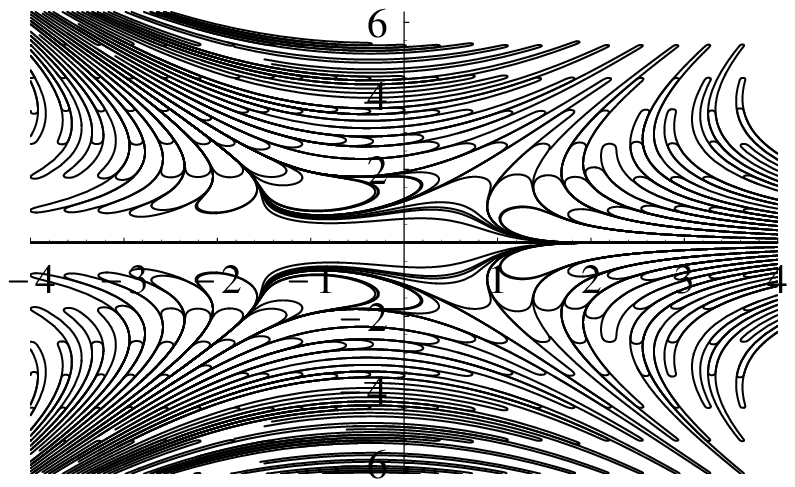)\par

\break

\head Accumulation of accumulating points
\endhead

\subsubhead Figure 7. $\Phi=e^{e^z} dz^2$ (See Theorem~4.2)
\endsubsubhead
Finally, this is an example about accumulation of
accumulations.  The image has infinitely many accumulation
points marked by dots outside the unit circle,
which in turns accumulate at $-1$.

\centredpicture 10cm by 7cm (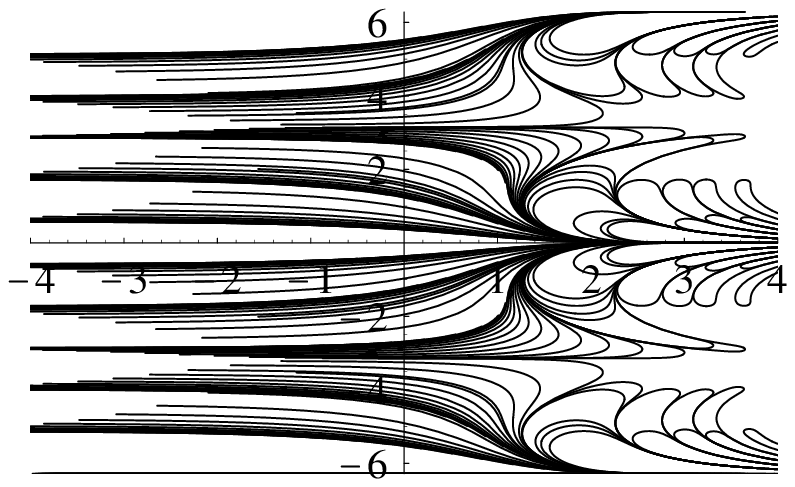)
\par\vskip 3ex plus 1ex minus 1ex\par
\centredpicture 7cm by 7cm (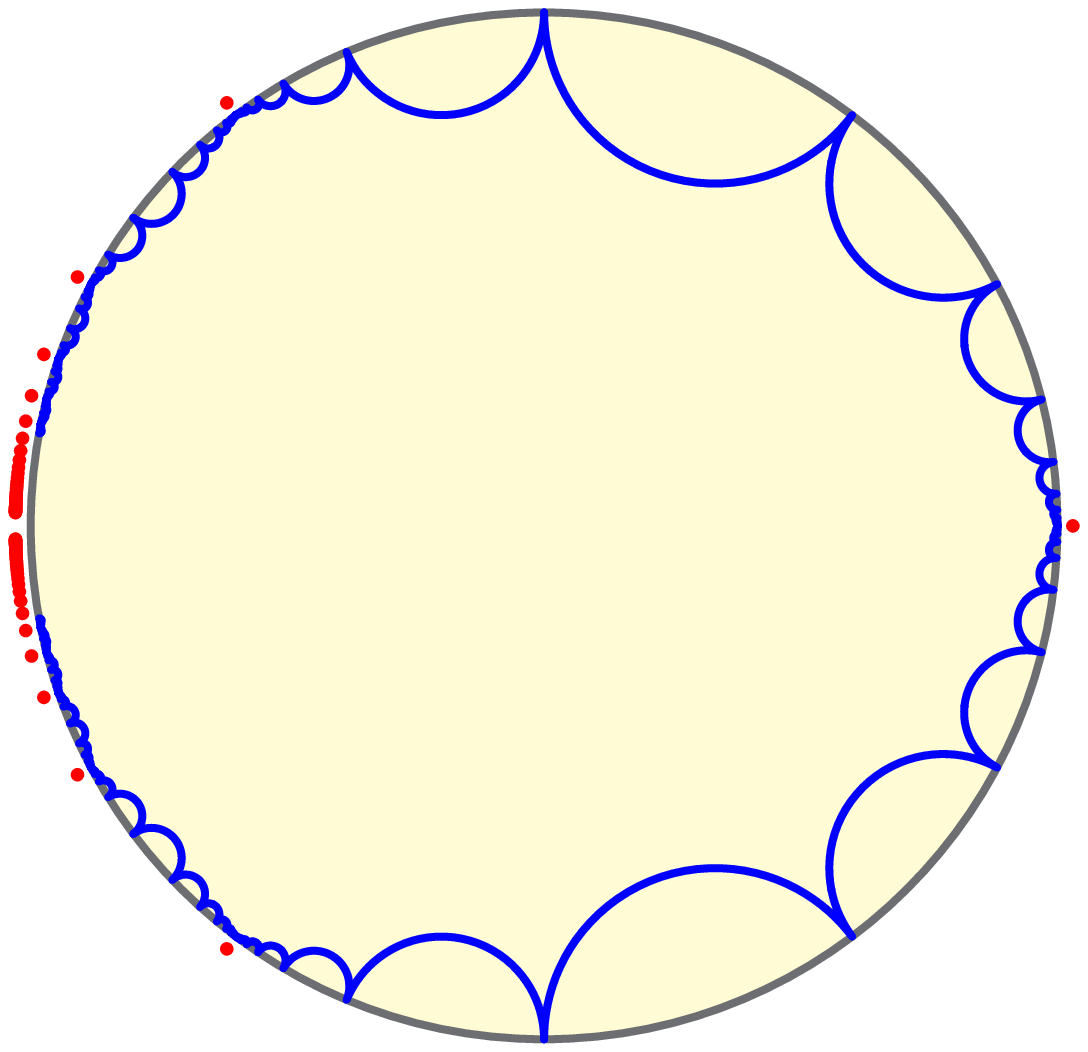)

\filbreak

\end